\documentclass[bj]{imsart}

\RequirePackage{amsthm,amsmath,amsfonts,amssymb}
\RequirePackage[numbers]{natbib}
\RequirePackage[colorlinks,citecolor=blue,urlcolor=blue]{hyperref}
\RequirePackage{graphicx}
\usepackage{dsfont}
\usepackage{subfigure}
\DeclareMathOperator*{\argmin}{arg\,min}

\startlocaldefs
\theoremstyle{plain}

\newtheorem{theorem}{Theorem}[section]
\newtheorem{lemma}[theorem]{Lemma}
\newtheorem{corollary}{Corollary}[section]
\theoremstyle{remark}

\newtheorem{remark}{Remark}

\newcommand{\bw}{\mathbf{w}}
\newcommand{\bv}{\mathbf{v}}

\newcommand{\balpha}{\mathbf{\alpha}}
\newcommand{\bbeta}{\mathbf{\beta}}
\newcommand{\bgamma}{\mathbf{\gamma}}

\newcommand{\D}{{\mathcal{D}}}
\newcommand{\F}{{\mathcal{F}}}

\newcommand{\N}{\mathbb{N}}

\newcommand{\R}{\mathbb{R}}
\newcommand{\Z}{\mathbb{Z}}
\newcommand{\Rd}{\mathbb{R}^d}

\newcommand{\beq}{\begin{eqnarray*}}

\newcommand{\eeq}{\end{eqnarray*}}

\newcommand{\beqm}{\begin{eqnarray}}

\newcommand{\eeqm}{\end{eqnarray}}

\newcommand{\EXP}{{\mathbf E}}
\newcommand{\PROB}{{\mathbf P}}

\renewcommand{\bf}{\normalfont \bfseries}
\renewcommand{\it}{\normalfont \itshape}

\endlocaldefs

\begin{document}

\begin{frontmatter}
\title{Convergence rates for shallow neural networks learned by gradient descent}
\runtitle{Shallow neural networks learned by gradient descent}

\begin{aug}
\author[A]{\fnms{Alina} \snm{Braun}\ead[label=e1,mark]{braun@mathematik.tu-darmstadt.de}},
\author[A]{\fnms{Michael} \snm{Kohler}\ead[label=e2,mark]{kohler@mathematik.tu-darmstadt.de}},
\author[B]{\fnms{Sophie} \snm{Langer}\ead[label=e3]{s.langer@utwente.nl}}
\and
\author[C]{\fnms{Harro} \snm{Walk}\ead[label=e4]{walk@mathematik.uni-stuttgart.de}}
\address[A]{Department of Mathematics,
Technical University of Darmstadt,
\printead{e1,e2}}

\address[B]{Department of Applied Mathematics,
University of Twente,
\printead{e3}}

\address[C]{Department of Mathematics,
University of Stuttgart,
\printead{e4}}
\end{aug}

\begin{abstract}
 In this paper we analyze the $L_2$
  error of neural network regression estimates
  with one hidden layer. Under the assumption
  that the  Fourier transform of the regression function
  decays suitably fast, we show that an estimate, where
  all initial weights are chosen according to proper
  uniform distributions and where the weights are learned
  by gradient descent, achieves
  a rate of convergence of $1/\sqrt{n}$
    (up to a logarithmic factor).
  Our statistical analysis implies that the key aspect
  behind this result is the proper choice of the initial
  inner weights and the adjustment of the outer weights
  via gradient descent. This indicates that we can also
  simply use linear least squares to choose the outer
  weights. We prove a corresponding theoretical result
  and compare our
new linear least squares neural network estimate
  with standard neural network estimates via simulated
  data. Our simulations show that our theoretical considerations
  lead to an estimate with an improved performance in many cases.
\end{abstract}



\end{frontmatter}

\section{Introduction}
\label{se1}

\subsection{Scope of this article}
\label{se1sub1}
Understanding the success of neural networks in practical applications (see, e.g., Krizhevsky, Sutskever, and Hinton (2012), Kim (2014), Wu et al. (2016)
or Silver et al. (2017)) is arguably one of the most important goals of machine learning theory today.  The problem has been studied in a statistical context, by analyzing empirical risk minimizers based on various classes of neural networks and under different assumptions on the target function (see, e.g.,  Schmidt-Hieber (2020), Kohler, Krzy\.zak and Langer (2019), Suzuki and Nitanda (2019) and the literature cited therein).  A complementary line of work (see, e.g.,  Choromanska et al. (2015),  Allen-Zhu and Li (2019), Ghorbani et al. (2019) and the literature cited therein) deals with the optimization procedure of the networks and analyzes the gradient descent routine and its variants.  Although both areas partially contribute to the theoretical understanding of deep learning, they each omit important parts in their analysis. In particular, they either work in an ideal setting without any optimization error or analyze an optimization procedure without any statistical setting.  With the goal in mind to bridge the gap between these two research areas, the aim of this work is to answer the following question:

\begin{center}
\textit{
Can we derive rate of convergence results for neural network estimators learned by gradient descent in a nonparametric regression setting?}
\end{center}
For simplicity, we restrict ourselves to the class of \textit{shallow} neural networks, i.e., neural networks with only one hidden layer and assume regression functions with suitable decaying Fourier transforms (see \eqref{inteq11}).

\subsection{Nonparametric regression}
\label{se1sub2}
We consider
a $\Rd \times \R$--valued random
vector $(X,Y)$, where $X$ is the so--called observation vector
and $Y$ is the  so-called response. Assume the condition
$\EXP \{Y^2\}<\infty$.  We are interested in the functional correlation
between the response $Y$ and the observation vector $X$. Particulary, 
we are searching for a function $f^*: \mathbb{R}^d \to \R$ such that
\begin{align*}
\EXP\left\{|f^*(X) - Y|^2\right\} = \min_{f: \R^d \to \R} \EXP\left\{|f(X) - Y|^2\right\}.
\end{align*}
This minimum holds for $f^*(x) = m(x) = \EXP\{Y|X=x\}$ (see Section 1.1 in Gy\"orfi et al. (2002)),
which is why $m$ is the so-called regression function.
But, in applications the distribution of $(X,Y)$ is unknown.  A basic problem in statistics is to recover the 
unknown regression function $m$ from a sample of $(X,Y)$, i.e.,  a data set
\begin{equation}
  \label{inteq1}
\D_n = \left\{
(X_1,Y_1), \ldots, (X_n,Y_n) 
\right\},
\end{equation}
where
$(X,Y)$, $(X_1,Y_1)$, \ldots, $(X_n,Y_n)$ are independent and identically distributed (i.i.d.).
Particulary,  we are searching 
for an estimator
\[
m_n(\cdot)=m_n(\cdot, \D_n):\Rd \rightarrow \R
\]
of $m$ such that the so--called $L_2$ error
\[
\int |m_n(x)-m(x)|^2 {\PROB}_X (dx)
\]
is ``small'' (cf., e.g., Gy\"orfi et al. (2002)
for a systematic introduction to nonparametric regression and
a motivation for the $L_2$ error). 

\subsection{Least squares neural network estimators}
\label{se1sub3}
Neural networks try to mimic the human brain in order to define
classes of functions. The starting point is a very simple model
of a nerve cell, in which some kind of thresholding is applied to
a linear combination of the outputs of other nerve cells. This leads
to functions of the form
\[
f(x)= \sigma \left(
\sum_{j=1}^d w_j \cdot x^{(j)} + w_0
\right)
\quad (x = (x^{(1)}, \dots, x^{(d)})^T \in \R^d),
\]
where we call $w_0$, \dots, $w_d \in \R$ the weights of the neuron
and $\sigma : \R \rightarrow \R$ the activation
function.
Traditionally, so--called squashing functions are chosen as activation
functions, which are nondecreasing
and satisfy $\lim_{x \rightarrow - \infty} \sigma(x)=0$
and
$\lim_{x \rightarrow  \infty} \sigma(x)=1$.
An example is
the so-called sigmoidal or logistic squasher
\begin{equation}
  \label{inteq2}
\sigma(x)=\frac{1}{1+\exp(-x)} \quad (x \in \R).
\end{equation}
Recently, also unbounded activation functions are used, e.g., the
ReLU activation function 
\begin{align*}
\sigma(x)=\max\{x,0\}.
\end{align*}
Some works like Sonoda, Ishikawa and Ikeda (2021), Sitzmann et al. (2020) and the literature cited therein,  consider periodic activation functions of the form $\sigma(t)=\sigma(t+T)$, where $T$ denotes the length of the period.
\\
\\
The most simple form of neural networks are shallow networks, i.e., neural networks
with one hidden layer, in which a simple linear combination
of the above neurons is used to define a function $f:\Rd \rightarrow \R$
by
\begin{equation}
  \label{inteq3}
f(x)= \sum_{k=1}^K \alpha_k \cdot
\sigma \left(
\sum_{j=1}^d \beta_{k,j} \cdot x^{(j)} + \beta_{k,0}
\right)
+ \alpha_0.
\end{equation}
Here $K \in \N$ is the number of neurons. The weights
$\alpha_k \in \R$ $(k \in \{0, \dots, K\})$, $\beta_{k,j} \in \R$
$(k\in \{1, \dots, K\}, j \in \{0, \dots, d\})$ are then fitted to the data
(\ref{inteq1}) in order to define an estimate of the
regression function. This can be achieved for example by 
applying the principle of least squares, i.e., by
defining the regression estimator $m_n$ by
\begin{equation}
  \label{inteq4}
  m_n(\cdot)
  =
  \argmin_{f \in \F}
  \frac{1}{n} \sum_{i=1}^n |Y_i - f(X_i)|^2,
\end{equation}
where $\F$ is the set of all functions of the form
(\ref{inteq3}) with a fixed number of neurons $K$ and fixed
activation function $\sigma$.
\\
\\
The rate of convergence
of shallow neural network regression estimates
has been analyzed in
Barron (1994) and McCaffrey and Gallant (1994).
Barron (1994) proved a dimensionless rate of $n^{-1/2}$
(up to some logarithmic factor), provided the Fourier transform
of the regression function
has a finite first
moment, which basically
requires that the function becomes smoother with increasing
dimension $d$ of $X$.
McCaffrey and Gallant (1994) showed for any $\epsilon >0$ a rate of $n^{-\frac{2p}{2p+d+5}+\varepsilon}$
in case of a $p$-times continuously differentiable regression function, but their study was restricted to the use of a certain cosine squasher as activation function. As to related work, we mention K\r{u}rková und Sanguinetti (2008) with further references.
\\
\\
In deep learning, neural networks with
several hidden layers are used to define classes of functions.
Here, the neurons are arranged in $\mathcal{L} \in \N$ layers,
where the $k_s \in \N$ neurons in layer $s \in \{2, \dots, \mathcal{L}\}$
get the output of the $k_{s-1}$ neurons in layer $s-1$ as
input, and where the neurons in the first layer are applied
to the $d$ components
of the input. We denote the weight between neuron $j$ in layer
$s-1$ and neuron $i$ in layer $s$ by $w_{i,j}^{(s)}$. This leads
  to the following recursive definition of a neural network
  with $\mathcal{L}$ layers and $k_s$ neurons in layer $s \in \{1, \dots, \mathcal{L}\}$:
  \begin{equation}
    \label{inteq5}
f(x) = \sum_{i=1}^{k_{\mathcal{L}}} w_{1,i}^{(\mathcal{L})}f_i^{(\mathcal{L})}(x) + w_{1,0}^{(\mathcal{L})}
\end{equation}
for some $w_{1,0}^{(\mathcal{L})}, \ldots, w_{1,k_L}^{(\mathcal{L})} \in \mathbb{R}$ and for $f_i^{(\mathcal{L})}$'s recursively defined by
\begin{equation}
  \label{inteq6}
f_i^{(s)}(x) = \sigma\left(\sum_{j=1}^{k_{s-1}} w_{i,j}^{(s-1)} f_j^{(s-1)}(x) + w_{i,0}^{(s-1)} \right)
\end{equation}
for some $w_{i,0}^{(s-1)}, \dots, w_{i, k_{s-1}}^{(s-1)} \in \mathbb{R}$,
$s \in \{2, \dots, \mathcal{L}\}$,
and
\begin{equation}
  \label{inteq7}
f_i^{(1)}(x) = \sigma \left(\sum_{j=1}^d w_{i,j}^{(0)} x^{(j)} + w_{i,0}^{(0)} \right)
\end{equation}
for some $w_{i,0}^{(0)}, \dots, w_{i,d}^{(0)} \in \mathbb{R}$.
\\
\\
The rate of convergence of least squares estimates
based on multilayer neural networks 
has been analyzed in Kohler and Krzy\.zak (2017),
 Imaizumi and Fukumizu (2018),
Bauer and Kohler (2019),
Kohler, Krzy\.zak and Langer (2019),
Suzuki and Nitanda (2019),
Schmidt-Hieber (2020) and Kohler and Langer (2021).
One of the main results obtained in this context shows
that neural networks
can achieve some kind of dimension reduction,
provided the regression
function is a composition of (sums of)
functions, where the input dimension of each of the functions 
is at most $d^*<d$
(see Kohler and Langer (2020) for a motivation of such a function class).
In Kohler and Krzy\.zak (2017) it was shown that
in this case
 suitably defined least squares estimates based on multilayer
neural networks achieve the rate of convergence  $n^{-2p/(2p+d^*)}$
(up to some logarithmic factor) for $p \leq 1$.
This result also holds for $p>1$
provided the squashing function is suitably
smooth
as was shown in
Bauer and Kohler (2019). Schmidt-Hieber (2020) showed the surprising result
that this is also true for neural
networks which use the non-smooth ReLU activation function.
In Kohler and Langer (2021) it was shown that
such results also hold for very simply
constructed fully connected feedforward neural networks.
Kohler, Krzy\.zak and Langer (2019)
considered regression functions with low local dimensionality 
and demonstrated that neural networks
are also able to circumvent  the curse of dimensionality
in this context.
Results regarding the estimation of regression functions
which are piecewise polynomials
having partitions with rather general smooth boundaries
by neural networks
 have been
 derived in Imaizumi and Fukumizu (2018).
 That neural networks can also
achieve a dimension reduction in Besov spaces
was shown in
Suzuki and Nitanda (2019).

\subsection{Gradient descent}
\label{se1sub4}
In Subsection \ref{se1sub3} the neural network regression estimates
are defined
as  nonlinear least squares estimates, i.e., as functions which minimize
the empirical $L_2$ risk
over nonlinear classes of neural
networks. In practice, it is usually not possible to find this
global minimum and 
one tries to find a local minimum using, for instance, the
gradient descent algorithm.
\\
\\
Denote by $f_{net,\bw}$ the neural network defined by
(\ref{inteq5})--(\ref{inteq7}) with weight vector
\[
\bw=(w_{j,k}^{(s)})_{s = 0, \dots, \mathcal{L}, j=1, \dots, k_{s+1}, k=0, \dots, k_{s}}
\]
(where we set $k_0=d$ and $k_{\mathcal{L}+1}=1$), and set
\begin{equation}
  \label{inteq8}
F(\bw)
=
\frac{1}{n} \sum_{i=1}^n
|Y_i - f_{net, \bw}(X_i)|^2.
\end{equation}
Now gradient descent is used to minimize
(\ref{inteq8}) with respect to $\bw$.
 Here, set
\begin{equation}
  \label{se2eq2}
\bw(0)=\bv
\end{equation}
for some  (usually randomly chosen)
initial weight vector $\bv$
and define
\begin{equation}
  \label{inteq9}
\bw(t+1) = \bw(t) - \lambda_n \cdot \nabla_\bw F (\bw(t))
\end{equation}
for
$t \in \{0,1, \dots, t_{n}-1\}$, where $\lambda_n>0$ is the stepsize
and $t_n \in \N$ is the number of performed gradient descent steps.
The  estimate is then defined by
\begin{equation}
  \label{inteq10}
  {m}_n(\cdot)
  =
  f_{net,\bw(t_n)}(\cdot).
\end{equation}

\subsection{Main results}
\label{se1sub5}
The main results in this article are threefold:
Firstly, we analyze the rate of convergence of a shallow neural
network regression estimate, where
the weights are learned by gradient descent. Here we assume
that the Fourier transform
\begin{align*}
\hat{F}: \Rd \rightarrow \mathbb{C}, \quad
\hat{F}(\omega)
=
\frac{1}{(2 \pi)^{d/2}} \cdot
\int_{\Rd}
\exp(- i \cdot \omega^T  x) \cdot m(x) \, dx
\end{align*}
of the regression function $m(x)=\EXP\{Y|X=x\}$
satisfies
\begin{equation}
  \label{inteq11}
 |\hat{F}(\omega)| \leq
    \frac{c_1}{\| \omega \|^{d+1} \cdot (\log \| \omega\|)^2 }
    \quad
    \mbox{for all } \omega \in \Rd \mbox{ with }
    \| \omega \| \geq 2
\end{equation}
for some $c_1 \in \R_+$.
We show that if we use the logistic squasher as
the activation function, if we choose the initial weights
of the neural network
randomly
from some proper uniform distributions,
and if we perform (up to some logarithmic factor) $n^{1.75}$
gradient descent steps with step size of order
$1/n^{1.25}$ (up to some logarithmic factor)
applied to some properly regularized empirical $L_2$ risk,
then
a truncated version of the estimate achieves (up to
some logarithmic factor) the rate of convergence
$1/\sqrt{n}$.
This shows that the classical result from Barron (1994)
also holds for a neural network estimate learned
by gradient descent. Surprisingly, in this result
a single random initialization of the weights is sufficient.
Furthermore our proof clarifies that this result mainly holds
because of the proper initialization of the weights and
because of the good adjustment of the outer weights
of the neural network during gradient descent.
We also establish a
minimax lower bound for the rate of convergence.
It reveals that for large $d$
the obtained rate of convergence is in its exponent $-1/2$
close to the optimal minimax rate of convergence. 
\\
Secondly,  we use our theoretical findings 
 to simplify
our estimate. Due to the fact that the optimization of
the inner weights by gradient descent is not necessary
in our result, it is evident that it should suffice
to minimize the outer weights of the neural network.
But this is (for fixed inner weights), in fact, a linear
least squares problem, for which the optimal weights
can easily be computed by solving a linear equation system.
We define a corresponding linear least squares estimator
with randomly selected inner weights, and show that for this
estimator the same rate of convergence result holds as for our
neural network estimator based on gradient descent.
The big advantage of this estimator is that it can be computed
much faster in applications. 
\\
Thirdly,  we compare our (theoretically motivated) estimator to classical 
shallow neural networks learned
by gradient descent on simulated data. In many cases we see a clear 
outperformance of our estimator over the classical ones.  

%
%
%

\subsection{Discussion of related results}
\label{se1sub6}
Our result shows that it is possible to extend
the classical result from Barron (1994) to the case
of a neural network estimator learned by gradient descent.
In contrast to Barron (1994), in which it was assumed
that the Fourier transform $\hat{F}$
of the regression function
has a finite first moment, i.e., 
\begin{align}
\label{Barronclass}
\int_{\Rd} \| \omega\| \cdot |\hat{F}(\omega)| \, d \omega < \infty,
\end{align}
we need the slightly stronger assumption
(\ref{inteq11}). 
\subsubsection{On related proof strategies}
Stone (1982) showed that the optimal minimax rate of convergence
for estimation of a $p$--times continuously differentiable
regression function is $n^{-2p/(2p+d)}$.  For fixed $p$ and 
increasing dimension $d$ this optimal rate gets worse in high
dimensions (so--called curse of dimensionality).
%
The rate $1/\sqrt{n}$ derived by Barron (1994) and also
in this paper is independent of the dimension and does
consequently not suffer from the curse of dimensionality. This
is due to the fact that the existence of a first moment
of the Fourier transform of the regression function
basically requires that
the smoothness of the regression function increases in case
of a growing dimension (cf., Remarks \ref{re3} and \ref{re4} below).
\\
For his statistical investigation, Barron (1994) used a result of Barron (1993) on the rate of approximation of a function with finite first moment of its Fourier transform by a shallow neural network. Barron (1993) obtained his deterministic approximation result by a probabilistic argument of Maurey (see Pisier (1980)).  This approach was analyzed and modified by Igelnik and Fao (1995) and motivated them to propagate shallow neural networks with random $d$-dimensional weight vectors $\beta_j$ and biases $\gamma_j$ ($j\in \{1, \dots, K_n\}$). Here the non-linear optimization problem in Barron (1994) is reduced to a quadratic optimization problem on the outer weights $\alpha_k$ ($k \in \{1, \dots, K_n\}$). As an approximation result, the authors established rate of convergence of the mean squared error for Lipschitz continuous functions (see also Huang et al. (2006)). \\
Section 3 of the present paper deals with the investigation of statistical learning of such a neural network, combining methods of empirical process theory and stochastic type approximation, modifying and partially weakening Barron's (1994) first moment condition.  Beside the estimation and approximation error an optimization error is taken into account, i.e., networks trained by gradient descent are considered. In a rather general framework Rahimi and Recht (2006) obtained results on learning random neural networks guaranteeing assertion validity with high probability.
\\
Under sharpened Barron conditions (higher moment conditions), which is satisfied, among other things, by solution functions of Kolmogorov PDEs, Goron (2021) established approximation, estimation and optimization error bounds for shallow neural network estimators with ReLU activation function and randomly generated internal weights and biases (see also Goron et al. (2020) with further literature).
\\
In a more practically oriented result, Dudek (2019) proposed a method for shallow neural network regression estimates on how to choose the range of random inner weights and biases depending on the input data and the shape of the activation functions.  The main result of our paper concerns a well-defined size of the range of the inner weights and biases, while in the proof, particulary in application of Lemma 5.1, the special shape of the logistic squasher is taken into account. It should be mentioned that for multilayer neural networks with randomly chosen inner weights and biases, Widrow et al. (2013) presented a gradient descent method for determing the outer weights.

\subsubsection{On the optimization error of neural networks}
There exist quite a few papers which try to show that neural network
estimators learned by gradient descent have nice theoretical properties.
The most popular approach in this context is the so--called
landscape approach. 
Choromanska et al. (2015)
used random matrix theory to
derive a heuristic argument showing
that the  risk  of  most  of  the  local  minima  of the
empirical $L_2$ risk is  not  much
larger  than  the  risk  of  the  global  minimum.  For networks with 
linear or quadratic activation function this claim could be 
validated, see, e.g. , Arora et al. (2018), 
Kawaguchi (2016),
and Du and Lee (2018).
However, these networks do not have good approximation properties.
Consequently, it is not possible to derive comparable convergence rates from these 
results as in our work.
Du et al. (2018) 
analyzed gradient descent applied to shallow neural networks in case of a Gaussian input distribution. 
But they used the expected gradient instead of the true gradient 
in their gradient descent routine
and therefore 
their result cannot be applied to derive the same convergence rates as in our work.
Liang et al. (2018)
applied gradient descent to a modified loss function in classification, 
where it is assumed that the data can be interpolated by a neural network.
Here, the second assumption is not satisfied in nonparametric regression
and it is unclear whether the main idea (of simplifying the estimation by a modification of the loss function) 
can also be used in a regression setting.  Brutzkus et al.  (2018) prove that two-layer networks with ReLU activation
function can learn linearly-separable data using stochastic gradient descent.  Andoni et al. (2014) also consider 
two-layer networks and analyze the sample complexity of these networks for learning multidimensional
 polynomial functions of finite degree.  But this result is based on exponential activation functions.
 For an overview of the literature concerning neural networks learned by gradient descent we also refer
 to Poggio, Banburski and Liao (2020).
\\
Our result can be understood as a confirmation of the conjecture in the
landscape approach in case of shallow neural networks.
We show that with our random initialization of the inner
weights of the neural network,
with high probability they are chosen such that there exist
values for the outer weights such that the corresponding neural
network has a small empirical $L_2$ risk. So, if we define the local
minima of the empirical $L_2$ risk as the minima which we get
if we just choose the outer weights optimally and keep the
values of the inner weights, then indeed most of the local minima
of the empirical $L_2$ risk have a small value.  This
is related to
 the assertion of Goodfellow,  Bengio and Courville (2015, pp. 3-5), 
 who mention that machine learning algorithms
 heavily depend on the representation of the data.
 In particular, they consider the ability of deep learning to learn a good hierarchical 
 representation of the data as a key aspect of its success.
In our result the inner representation of the data used by our network
depends on the randomly chosen inner weights, which are applied to
the activation function. 
 Hence, in our result the key feature of the
 neural networks is \textit{representation guessing} instead of
 \textit{representation learning}.
\\
 For a related topic, i.e., estimation of regression functions
 by generalizations of two layer
 radial basis function networks, the asymptotic
 behaviour of the gradient descent was analyzed in
 Javanmard, Mondelli and Montanari (2021) by using a so--called
 Wasserstein gradient descent approach. It remains unclear whether
the approach can be extended to classical neural
 networks. In particular,  it is unclear whether the results about shallow networks as in our article
 can be derived by this approach.
\\
\subsubsection{On results of overparametrized neural networks}
 Recently it was shown in quite a few papers that
 in case of  suitably overparameterized neural networks
 gradient descent can find the global minimum of the empirical $L_2$ risk,
 cf., e.g.,
Kawaguchi and Huang (2019),
Allen-Zhu, Li and Song (2019),
Allen-Zhu, Li and Liang (2019), 
Arora et al. (2019a, 2019b), 
Du et al. (2018), 
Li and Liang (2018) and
Zou et al. (2018).
However, Kohler and Krzy\.zak (2019)
presented
a counterexample demonstrating
that overparameterized neural networks,
which basically interpolate the training data,
in general do not generalize well.
In this counterexample the regression function is constant zero
and hence satisfies the assumption on the regression function
imposed in our paper.
In particular, this shows that
results similar to the ones in our paper
cannot be concluded from the papers cited above.
We would also like to stress that our estimator does not
use an overparameterization, because the numbers of weights
of our neural networks in the theorems below are
much smaller than the sample size.
\\
Another approach to analyze overparameterized neural networks
is the the so--called kernel approach
(cf. Jacot, Gabriel and Hongler (2020)
and the literature cited in Woodworth et al. (2020)).
Here, neural networks are approximately described
by kernel methods and a gradient descent in continuous
time modelled by a differential equation leads to
the so--called neural tangent kernel, which depends
on the time. The asymptotic behaviour of this
neural tangent kernel has been analyzed in
Jacot, Gabriel and Hongler (2020), which leads to
an asymptotic approximation of neural networks.
Unfortunately this asymptotic approximation does not
imply how the finite neural networks behave during
learning.

\subsubsection{On the generalization error of neural networks}
The generalization of neural networks can also be analyzed
within the classical Vapnik Chervonenkis theory (cf., e.g.,
Chapters 9 and 17 in  in Gy\"orfi et al. (2002)).
Here, the complexity of the underlying function spaces
is measured by covering numbers, which can be bounded using
the so--called Vapnik-Chervonenkis dimension (cf., e.g., Bartlett
et al.
(2019)).
However,
the resulting upper bounds on the generalization
error might be too rough as during gradient descent
the neural network estimator does not necessarly attend all functions
from the underlying function space.  One might sharpens the bound
by using the 
so--called Rademacher complexity (cf., Koltchinski (2004)).
For networks with quadratic activation function
this has already been successfully done in Du and Le (2018),
but unfortunately such neural networks do not have good approximation properties and similar results as in our work can therefore certainly not be derived.
Also, we would like
to stress that in our result we indeed analyze the
generalization of neural networks within the classical
Vapnik Chervonenkis theory.
\subsection{Notation}
\label{se1sub7}
Throughout the paper, the following notation is used:
The sets of natural numbers, natural numbers including $0$, real numbers,
nonegative real numbers and complex numbers
are denoted by $\N$, $\N_0$, $\R$, $\R_+$ and $\mathbb{C}$, respectively. For $z \in \R$, we denote
the smallest integer greater than or equal to $z$ by
$\lceil z \rceil$ and the largest integer smaller or equal to $z$ by 
$\lfloor z \rfloor$. 
Let $D \subseteq \R^d$ and let $f:\R^d \rightarrow \R$ be a real-valued
function defined on $\R^d$.
We write $x = \arg \min_{z \in D} f(z)$ if
$\min_{z \in \D} f(z)$ exists and if
$x$ satisfies
$x \in D$ and $f(x) = \min_{z \in \D} f(z)$.
The Euclidean norm of $x \in \Rd$
is denoted by $\|x\|$.
For $f:\R^d \rightarrow \R$
\[
\|f\|_\infty = \sup_{x \in \R^d} |f(x)|
\]
is its supremum norm.
$S_r$ denotes the ball with radius $r$ in $\Rd$
and center $0$
(with respect to the Euclidean norm).
We define the truncation operator $T_{\kappa}$ with level $\kappa > 0$ as
\begin{equation*}
T_{\kappa}u =
\begin{cases}
u \quad &\text{if} \quad |u| \leq\kappa\\
\kappa \cdot {\rm sign}(u) \quad &\text{otherwise}.
\end{cases}
\end{equation*}
Constants are designated and numbered $c_1, c_2, \dots$.   Each constant is assumed to be non-negative and, unless otherwise stated, absolute.
\subsection{Outline}
\label{se1sub8}
In Section \ref{se2} we present our main result concerning
the rate of convergence of a shallow neural network estimator
 learned by gradient descent. In Section \ref{se3}
we show that the same rate of convergence can also be achieved
by a  linear least squares estimator with much simpler
computation.
In Section \ref{se4} we compare the finite sample size behaviour
of our linear least squares estimate via simulated data. Section \ref{se5} contains
the proof of a key auxiliary result concerning the approximation
error of shallow neural networks with randomly chosen inner weights,
and the outline of the proofs of Theorem \ref{th1} and
Theorem \ref{th2}.
The complete proofs of our main results are given in the supplement.

\section{A neural network estimate learned by gradient descent}
\label{se2}
In this sequel we analyze shallow neural networks with 
$K_n$ hidden neurons and a constant term. As activation function 
we choose the logistic
squasher (\ref{inteq2}).
The networks are defined by
\begin{equation}
  \label{se2eq6}
  f_{net,\bw}(x)
  =
  w_{1,0}^{(1)}+
  \sum_{j=1}^{K_n}
  w_{1,j}^{(1)} \cdot \sigma\left(
  \sum_{k=1}^d w_{j,k}^{(0)} \cdot x^{(k)}
  + w_{j,0}^{(0)}\right)
  =
  \alpha_0
  +
  \sum_{j=1}^{K_n}
  \alpha_j \cdot \sigma\left( \beta_j^T \cdot x + \gamma_j\right)
\end{equation}
where $K_n \in \N$,
$\alpha_i, \gamma_i \in \R$, 
$\beta_i=(\beta_{i,1}, \dots, \beta_{i,d})^T \in \Rd$ $(i \in \{1, \dots, K_n\})$
and
\[
\bw=(w_{j,k}^{(l)})_{j,k,l}=(\balpha,\bbeta,\bgamma)
=(\alpha_0,\alpha_1, \dots, \alpha_{K_n},\beta_1,\dots,\beta_{K_n},
\gamma_1,\dots,\gamma_{K_n})
\]
is the vector of the $D_n=1+K_n \cdot (d+2)$ 
weights of the neural network $f_{net,\bw}$.
\\
\\
We learn the weight vector  by minimizing the regularized
least squares criterion
\begin{equation}
  \label{se2eq1}
F(\bw)
=
\frac{1}{n} \sum_{i=1}^n
|Y_i - f_{net, \bw}(X_i)|^2
+
\frac{c_2}{K_n} \cdot \sum_{k=0}^{K_n} \alpha_k^2,
\end{equation}
where $c_2>0$ is an arbitrary constant.
Minimization of (\ref{se2eq1}) with respect to $\bw$
is a nonlinear least squares problem, for which we use
gradient descent. Here we set
\begin{equation}
  \label{se2eq2}
\bw(0)=\bv,
\end{equation}
where the  
initial weight vector
\[
\bv=
(\alpha_0^{(0)},\alpha_1^{(0)}, \dots, \alpha_{K_n}^{(0)},
\beta_1^{(0)},\dots,\beta_{K_n}^{(0)},
\gamma_1^{(0)},\dots,\gamma_{K_n}^{(0)})
\]
is chosen such that
\begin{equation}
  \label{se2eq2b}
  |\alpha_0^{(0)}| \leq c_3
  \quad \mbox{and} \quad
  |\alpha_k^{(0)}| \leq \frac{c_4}{K_n}
  \quad (k \in \{1, \dots, K_n\})
\end{equation}
holds for constants $c_3, c_4 > 0$ and such that
$\beta_1^{(0)},\dots,\beta_{K_n}^{(0)},
\gamma_1^{(0)},\dots,\gamma_{K_n}^{(0)}$
are independently distributed with
$\beta_1^{(0)},\dots,\beta_{K_n}^{(0)}$ uniformly
distributed on $\{ x \in \Rd \, : \, \|x\|=B_n\}$
and
$\gamma_1^{(0)},\dots,\gamma_{K_n}^{(0)}$ uniformly distributed
on  $[-B_n \cdot \sqrt{d} ,B_n \cdot \sqrt{d}]$
(where $B_n>0$ is defined in Theorem \ref{th1} below).
We define
\begin{equation}
  \label{se2eq3}
\bw(t+1) = \bw(t) - \lambda_n \cdot \nabla_\bw F (\bw(t))
\end{equation}
for
$t \in \{0,1, \dots, t_{n}-1\}$. Here, $\lambda_n>0$ is the stepsize
and $t_n \in \N$ is the number of performed gradient descent steps. 
Both are defined in Theorem \ref{th1}
below.
Our estimator is then given by
\begin{equation}
  \label{se2eq4}
  \tilde{m}_n(\cdot)
  =
  f_{net,\bw(t_n)}(\cdot)
\end{equation}
and
\begin{equation}
  \label{se2eq5}
  m_n(x)=T_{\kappa_n}  \tilde{m}_n(x)
  \end{equation}
where $\kappa_n= c_5 \cdot \log n$.
The truncation operator is necessary for theoretical reasons. Later we apply results from empirical process theory to bound the covering number and the VC dimension of our function space of shallow neural networks.  Here boundedness of the function space is needed.  As an alternative one could directly restrict the class of shallow neural networks by imposing a sup-norm bound on all functions in the space (see, e.g., Schmidt-Hieber (2020)).  In case of restricted activation functions, like sigmoid or tangens hyperbolicus, one could also impose restrictions on the outer weights. But as we are applying gradient descent,  this would mean that we would have to check the weights after each gradient step.

\begin{theorem}
  \label{th1}
  Let $(X,Y)$ be an $[0,1]^d \times \R$--valued random vector such that
  \begin{equation}
    \label{th1eq1}
    \EXP \left\{ \exp(c_6 \cdot Y^2) \right\} < \infty
  \end{equation}
  holds for some constant $c_6>0$ and assume that
the
  corresponding regression function
  $m(x)=\EXP\{Y|X=x\}$  is bounded, satisfies
  \[
\int_{\Rd} | m(x)| \, dx < \infty,
  \]
  and that its Fourier transform $\hat{F}$  satisfies
 \begin{equation}
   \label{th1eq2}
    |\hat{F}(\omega)| \leq
    \frac{c_1}{\| \omega \|^{d+1} \cdot (\log \| \omega\|)^2 }
    \quad
    \mbox{for all } \omega \in \Rd \mbox{ with }
    \| \omega \| \geq 2
 \end{equation}
 for some $c_1>0$.
 Set
  \[
  K_n = \lceil c_7 \cdot \sqrt{n} \rceil, \quad
  B_n= \frac{1}{\sqrt{d}} \cdot (\log n)^2 \cdot K_n \cdot n^2,
  \]
  \[
   L_n 
   =
   c_8 \cdot (\log n)^6 \cdot K_n^{5/2}, \quad
  	\lambda_n = \frac{1}{L_n}  
  \]
  and 
  \[
  	t_n = \lceil K_n \cdot (\log n)^2 \cdot L_n \rceil,
        \]
        let $\sigma$ be the logistic squasher
  and define
  the estimator $m_n$ of $m$ as in \eqref{se2eq5}.
  Then one has for $n$ sufficiently large
    \begin{eqnarray*}
      &&
      \EXP \int | m_n(x)-m(x)|^2 \PROB_X (dx)
    \leq
    c_9 \cdot   
    (\log n)^{4} \cdot
    \frac{1}{\sqrt{n}}.
    \end{eqnarray*}
  \end{theorem}

\begin{remark}
  \label{re1}
    The computation of the estimator in Theorem \ref{th1}
    requires 
    \[
    t_n \leq c_{10} \cdot (\log n)^8 \cdot (\sqrt{n})^{7/2}
    =
    c_{10} \cdot (\log n)^8 \cdot n^{7/4} 
    \]
    many (i.e., up to a logarithmic factor only $n^{1.75}$ many)
    gradient descent steps and only one initialization
    of the starting weights.
\end{remark}

\begin{remark}
  \label{re2}
        Condition (\ref{se2eq2b}) is in particular
        satisfied if we set $\alpha_k^{(0)}=0$ $(k \in \{0, \dots, K_n\})$
        or if we choose  $\alpha_k^{(0)}$ $(k \in \{0, \dots, K_n\})$ 
         independently uniformly distributed on the interval
        $[-c_{11} /K_n,c_{11}/K_n]$ for some constant $c_{11} > 0$.
        \end{remark}
        
\begin{remark}
  \label{re3}
          Let $m:\Rd \to \R$ be $k$-times continously differentiable with Lebesgue integrable $k$-th partial derivatives. Then by standard Fourier analysis
          (cf., e.g., Proposition 4.5.3 in Epstein (2008))
        \begin{align*}
        |\hat{F}(\omega)| \leq \frac{c_{12}}{(1+\|\omega\|)^k}, \quad \omega \in \Rd.
        \end{align*}
        If $k=d+2$, then (\ref{th1eq2}) is fulfilled.
        This sufficient
condition is far from being necessary (see Theorem \ref{th3} below).
        \end{remark}

        \begin{corollary}
          \label{co1}
          If in Theorem \ref{th1} the regression function $m:\Rd \to \R$ is
radially symmetric and
          Lebesgue integrable, i.e., $m(x) = m^*(\|x\|)$, $x \in \Rd$, for some $m^*:\R_+ \to \R$ satisfying
        \begin{align*}
        \int_{\R_+} m^*(r)r^{d-1} dr < \infty, 
        \end{align*}
        then \eqref{th1eq2} can be replaced by
        \begin{align}
        \label{c1eq1}
        \int \|\omega\| \cdot |\hat{F}(\omega)| d\omega < \infty.
        \end{align}
        \end{corollary}
        
        \begin{remark} 
                Let $\partial S_r$ denote the surface of $S_r$ and $w_r = const(d) \cdot r^{d-1}$ denote its $(d-1)$- dimensional Lebesgue measure. Then Barron's (1994) condition \eqref{c1eq1} means 
                \begin{align*}
                \int_{\mathbb{R}_+} r \left(\int_{\partial S_r} |\hat{F}|d\sigma\right) dr < \infty,
                \end{align*}
            and \eqref{th1eq2} means
            \begin{align*}
            |\hat{F}(\omega)| \leq \frac{c_1}{r^{d+1} (\log r)^2} \ \text{for} \ \|\omega\|=r \geq 2r.
            \end{align*}
            Obviously, \eqref{th1eq2} implies \eqref{c1eq1}.
            If the Barron condition is sharpened to
            \begin{align*}
            \int_{\mathbb{R}_+} r w_r \sup_{\omega \in \partial S_r} |\hat{F}(\omega)| dr < \infty
            \end{align*}
            and $ \sup_{\omega \in \partial S_r} |\hat{F}(\omega)|$ is assumed non-increasing,
            then 
            \begin{align*}
            \sup_{\omega \in \partial S_r} |\hat{F}(\omega)| = o(r^{-d-1}) \quad (r \to \infty),
            \end{align*}
            which up to a logarithmic term corresponds to \eqref{th1eq2}. An analogous conclusion concerns $\int|\hat{F}(\omega)|d\omega < \infty$ and \eqref{re5eq1} in Section 3. Both conclusions immediately follow from the well-known fact (compare Olivier's theorem on infinite series) that for $q > 0$ and non-increasing $h: \mathbb{R}_+ \to \mathbb{R}_+$ finiteness of $\int_{\mathbb{R}_+} r^q h(r) dr$implies
            \begin{align*}
            h(r) = o(r^{-q-1}) \quad (r \to \infty).
            \end{align*}
            In the particular situation of Corollary \ref{co1}, $\hat{F}$ is radially symmetric (cf., Theorem 4.5.3 in Epstein (2008)), and especially $\hat{F}$ with
            \begin{align*}
            \hat{F}(\omega) = \frac{const}{r^{d+1} (\log r) (\log (1+\log r))^2} \ \text{for} \ \|\omega\|=r \geq 2
            \end{align*}
            satisfies  \eqref{c1eq1}, but not \eqref{th1eq2}.
            \end{remark}
            
        \begin{remark}
          \label{re4}
          Let $k \in \N$ and assume that
          $m: \Rd \to \R$ is  $k$-times continuously differentiable and
          that the $k$-th partial derivatives are square Lebesgue integrable.
          Noticing that the Fourier transform of the square integrable function
          \[
\frac{\partial^k m}{\partial (x^{(j)})^k} 
          \]
          is the function
          \[
\omega \mapsto ( -i \cdot \omega^{(j)})^k \cdot \hat{F} (\omega)
\]
(cf., e.g., Proposition 4.5.3 in Epstein (2008))
and thus, by Parseval's formula
(cf., e.g., Theorem 4.5.2 in Epstein (2008)
and Section VI.2 in Yosida (1968)),
\begin{equation}
  \label{re4eq1}
\int
\left(
\frac{\partial^k}{\partial (x^{(j)})^k } m(x)
\right)^2
dx
=
\int
|\omega^{(j)}|^{2k} \cdot |\hat{F} (\omega)|^2 \, d \omega,
\end{equation}
one obtains
        \begin{align*}
        \int \|\omega\|^{2k} \cdot |\hat{F}(\omega)|^2 d\omega < \infty.
        \end{align*}
        If $k=\lfloor \frac{d}{2} \rfloor +2$ (weaker than in Remark \ref{re3}), then \eqref{c1eq1} is fulfilled, because the Cauchy-Schwarz inequality yields
        \begin{align*}
        &\left(\int\|\omega\| \cdot |\hat{F}(\omega)|d\omega\right)^2\\
        \leq & \int \frac{1}{(1+\|\omega\|)^{2k-2}} d\omega \cdot \int (1+\|\omega\|)^{2k-2} \|\omega\|^2 |\hat{F}(\omega)|^2 d\omega\\
        \leq &
        c_{13} \cdot \int_0^\infty \frac{r^{d-1}}{(1+r)^{2k-2}}dr \cdot \int \|\omega\|^{2k} \cdot |\hat{F}(\omega)|^2 d\omega  <  \infty.
        \end{align*}
        This consideration can be found in
        Lee (1996), Chapter 7, pp. 69, 70.

        A simple example satisfying the above
assumption is  the radially symmetric function
$m:\Rd \to \R$ with
        \begin{align*}
        m(x) := \begin{cases}
        (1-\|x\|)^{k+1}, \quad \text{if} \ \|x\| \leq 1\\
        0, \quad \text{if} \ \|x\| > 1,
        \end{cases}
        \end{align*}
        where $k = \lfloor \frac{d}{2}\rfloor +2$.

        On the other side, condition (\ref{th1eq2}) of
        Theorem \ref{th1}
        implies that for same
$k=\lfloor \frac{d}{2} \rfloor +2$
  the
  right-hand side of (\ref{re4eq1}) is finite, the left-hand side of
  (\ref{re4eq1}) exists in a
distribution-theoretic sense - according to Yosida (1968), Section VI.2,
especially Corollary 1 for inverse of Fourier transform -
and (\ref{re4eq1}) holds.
        \end{remark}
        
        \begin{remark}
  It is an open problem whether one can extend Theorem \ref{th1} to other activation functions like ReLU or deeper network structures.  The most important trick in our result is that the internal weights change only slightly during the gradient descent.  This also follows from the properties of the sigmoid function which are no longer valid in the case of the unbounded ReLU function.  It is questionable whether, in the case of several hidden layers, all weights also change only slightly during the gradient descent.  Additionally one has to think about proper initalizations for all hidden layers.  A first step would be to analyze networks with two hidden layers.
        \end{remark}

        The rate derived in Theorem \ref{th1} is close to the
        optimal minimax rate of convergence as our next theorem
        about the lower bound shows.
\noindent       
      \begin{theorem}
      \label{th3}
      Let $c_1, c_6, c_{14}, c_{15} >0$ be sufficiently large, and
      let $\mathcal{D}$ be the class of all distributions, where
      \begin{itemize}
      \item[(1)] $X \in [0,1]^d$ a.s.
      \item[(2)] $\EXP\{\exp(c_6 \cdot Y^2)\} < \infty$
      \item[(3)] $m$ is bounded in absolute value by $c_{14}$
      \item[(4)] $\int_{\R^d} |m(x)| dx \leq c_{15}$
      \item[(5)]     $|\hat{F}(\omega)| \leq
    \frac{c_1}{\| \omega \|^{d+1} \cdot (\log \| \omega\|)^2 }$
    for all
    $\omega \in \Rd$
    with 
    $\| \omega \| \geq 2$
      \end{itemize}
      (where $\hat{F}$ is the Fourier transform of $m(\cdot)=\EXP\{Y|X=\cdot\})$. Then we have for $n$ sufficiently large 
      \begin{align*}
      \inf_{\hat{m}_n} \sup_{(X,Y) \in \mathcal{D}} \EXP\int|\hat{m}_n(x) - m(x)|^2 \PROB_X(dx) \geq c_{16} \cdot (\log n)^{-6} \cdot n^{-\frac{1}{2} - \frac{1}{d+1}}.
      \end{align*}
      \end{theorem}
      
      \begin{remark}
        In case that $d$ is large the exponent in the lower bound in Theorem \ref{th3}
        gets arbitrarily close to $-1/2$, i.e., the lower bound is close to the rate of convergence of Theorem \ref{th1}.
      \end{remark}
         
\section{A linear least squares neural network estimator}
\label{se3}
In this section we show that we can achieve the rate of convergence of Theorem \ref{th1} also by a simple linear least squares estimator,
where the underlying linear function space consists of shallow neural
networks with randomly chosen inner weights.

In order to define our function space, we start by choosing
$\beta_1,\dots,\beta_{K_n},
\gamma_1,\dots,\gamma_{K_n}$
independently distributed such that
$\beta_1,\dots,\beta_{K_n}$ are uniformly
distributed on $\{ x \in \Rd \, : \, \|x\|=B_n\}$
and
$\gamma_1,\dots,\gamma_{K_n}$ are uniformly distributed
on  $[-B_n \cdot \sqrt{d} ,B_n \cdot \sqrt{d}]$
(where $B_n>0$ is defined in Theorem \ref{th2} below).
Then we set
\begin{equation}
  \label{se3eq1}
  \F_n
  =
  \left\{
  f: \Rd \rightarrow \R \, : \,
  f(x)
  =
  \alpha_0
  +
  \sum_{j=1}^{K_n}
  \alpha_j \cdot \sigma( \beta_j^T \cdot x + \gamma_j)
\mbox{ for some } \alpha_0, \dots, \alpha_{K_n}  \in \R
  \right\}.
\end{equation}
Using this (random) linear function space we define our
estimator by
\begin{equation}
  \label{se3eq2}
  \tilde{m}_n(\cdot)
  =
  \argmin_{f \in \F_n}
  \frac{1}{n} \sum_{i=1}^n |Y_i - f(X_i)|^2
\end{equation}
and
\begin{equation}
  \label{se3eq3}
  m_n(x)=T_{\kappa_n}  \tilde{m}_n(x),
  \end{equation}
where $\kappa_n= c_5 \cdot \log n$.

\begin{theorem}
  \label{th2}
  Let $(X,Y)$ be an $[0,1]^d \times \R$--valued random vector such that
  (\ref{th1eq1})
  holds for some constant $c_6>0$ and assume that
the
  corresponding regression function
  $m(x)=\EXP\{Y|X=x\}$  is bounded, satisfies
  \[
\int_{\Rd} | m(x)| \, dx < \infty,
  \]
  and that its Fourier transform
$\hat{F}$ 
 satisfies (\ref{th1eq2})
  for some $c_1>0$.
 Set
  \[
  K_n = \lceil c_{17} \cdot \sqrt{n} \rceil  \quad \mbox{and} \quad
  B_n= \frac{1}{\sqrt{d}} \cdot (\log n)^2 \cdot K_n \cdot n^2,
  \]
  let $\sigma$ be the logistic squasher,
  choose
$\beta_1,\dots,\beta_{K_n},
\gamma_1,\dots,\gamma_{K_n}$
as above and define the estimator $m_n$ by
(\ref{se3eq1}), 
(\ref{se3eq2}) and 
(\ref{se3eq3}).
Then we have for $n$ sufficiently large
    \begin{eqnarray*}
      &&
      \EXP \int | m_n(x)-m(x)|^2 \PROB_X (dx)
    \leq
    c_{18} \cdot   
    (\log n)^{4} \cdot
    \frac{1}{\sqrt{n}}.
    \end{eqnarray*}
  \end{theorem}

\begin{remark}
  \label{re5}
  If we ignore logarithmic factors, then Theorem \ref{th3} implies
  that the rate of convergence in Theorem \ref{th2} is optimal
  up to the factor
  \[
n^{-\frac{1}{d+1}}.
\]
In applications $d$ is usually rather large, therefore, in our estimation, 
this factor has no practical relevance. However, from a mathematical
point of view it is rather unsatisfying if a regression estimator does
not achieve an optimal rate of convergence
at least up to some logarithmic factor.

We believe that with respect to the derived convergence rate $1/\sqrt{n}$, 
assumption (\ref{th1eq2}) is somewhat too strong meaning that our proof strategy (based on the result of Barron (1994))
seems to be not
suitable to derive an optimal rate of convergence.  The next corollary shows that 
slightly weakening \eqref{th1eq2} and modifying \eqref{se3eq3} leads to minimax optimal
rate of convergence result. In particular,  we set
\begin{eqnarray*}
  &&
  \F_n
  =
  \Bigg\{
  f: \Rd \rightarrow \R \, : \,
  f(x)
  =
  \alpha_0
  +
  \sum_{j=1}^{K_n}
  \alpha_j \cdot \sigma \left(
  B_n \cdot \left( Proj_{(-\pi,\pi]}( \beta_j^T \cdot x) + \gamma_j \right)
    \right)
    \\
    &&
    \hspace*{10cm}
\mbox{for some } \alpha_0, \dots, \alpha_{K_n}  \in \R
  \Bigg\},
\end{eqnarray*}
where the projection operator is defined by
\[
Proj_{(-\pi,\pi]}(z)=z+k \cdot 2 \cdot \pi
\]
and  $k=k(z) \in \Z$ is chosen such that
\[
z+k \cdot 2 \cdot \pi \in (-\pi,\pi]
\]
holds and
$\beta_1,\dots,\beta_{K_n},
\gamma_1,\dots,\gamma_{K_n}$ are
independently distributed such that
$\beta_1,\dots,\beta_{K_n}$ have the density
\[
\omega \mapsto     \frac{1}{4^{d+1}} \cdot 1_{\{\| \omega \| \leq 2\}} +
  \frac{c_{20}}{\|\omega\|^d \cdot (\log \|\omega\|)^2}    
 \cdot 1_{\{\| \omega \| > 2\}}
\]
with respect to the Lebesgue measure
(which for a proper choice of $c_{20}>0$ is indeed a density, if
(\ref{re5eq1}) holds) and
  $\gamma_1$,
\dots, $\gamma_{K_n}$ are uniformly distributed on $[-\pi,\pi]$.
The corresponding estimator $m_n$ is then defined as in \eqref{se3eq2} and \eqref{se3eq3} 
by 
\begin{equation}
\label{neqmnt}
  \tilde{m}_n(\cdot)
  =
  \argmin_{f \in \F_n}
  \frac{1}{n} \sum_{i=1}^n |Y_i - f(X_i)|^2
\end{equation}
and
\begin{equation}
\label{neqmn}
  m_n(x)=T_{\kappa_n}  \tilde{m}_n(x),
  \end{equation}
where $\kappa_n= c_5 \cdot \log n$.

\begin{corollary}
\label{cor31}
Assume that in Theorem \ref{th2} condition \eqref{th1eq2} is weakened to 
\begin{equation}
  \label{re5eq1}
    |\hat{F}(\omega)| \leq
    \frac{c_1}{\| \omega \|^{d} \cdot (\log \| \omega\|)^2 }
    \quad
    \mbox{for all } \omega \in \Rd \mbox{ with }
    \| \omega \| \geq 2
  \end{equation}
  for some $c_1>0$ and $m_n$ is defined as in \eqref{neqmnt} and \eqref{neqmn}.
Then 
\[
\sup_{(X,Y) \in \mathcal{D}} \EXP\int|m_n(x) - m(x)|^2 \PROB_X(dx) \leq c_{21} \cdot (\log n)^{6} \cdot \frac{1}{\sqrt{n}}.
      \]
\end{corollary}

\begin{remark}

According to Theorem \ref{th3} this rate is, up to a logarithmic factor, optimal. In particular, we have
  \begin{align*}
      \inf_{\hat{m}_n} \sup_{(X,Y) \in \mathcal{D}} \EXP\int|\hat{m}_n(x) - m(x)|^2 \PROB_X(dx) \geq c_{16} \cdot (\log n)^{-6} \cdot \frac{1}{\sqrt{n}}.
      \end{align*}
\end{remark}

\begin{remark}
At this point it should be emphasised that the weaker condition \eqref{re5eq1} leads to optimal rates but is no longer a subset of the Barron class , i.e., the class of all regression function where the Fourier transform
satisfies \eqref{Barronclass} (or \eqref{c1eq1}). This in turn means that we cannot consider the result as a special case of the Barron class, as is the case with the stronger condition.  
\end{remark}

\end{remark}

\section{Application to simulated data}
\label{se4}
This section provides a simulation-based comparison of our new linear least squares estimator with standard neural network estimators defined in the deep learning framework of Python's \textit{tensorflow} and \textit{keras}.  To implement our new estimator we compute in a first step the values of 
\begin{align*}
\sigma(\beta_j^T \cdot X_i + \gamma_j)
\end{align*}
for $i \in \{1, \dots, n_{learn}\}$, $j \in \{1, \dots, K_n\}$ and then solve a linear equation system for the values of $\alpha_0, \dots, \alpha_{K_n}$. In the initialization of the weights we choose the inner weights $\beta_1^{(0)}, \dots, \beta_{K_n}^{(0)}$ (according to the theoretical results) uniformly distributed on $\{x \in \Rd: \|x\|=B_n\}$ and the inner bias terms $\gamma_1^{(0)}, \dots, \gamma_{K_n}^{(0)}$ uniformly distributed on $[-B_n \cdot \sqrt{d}, B_n \cdot \sqrt{d}]$.  The values of $K_n$ and $B_n$ are chosen in a data-dependent way by splitting of the sample.   Here we use
$n_{learn}=\lceil \frac{4}{5} \cdot n \rceil$ realizations to train the estimator several times with different choices of $K_n$ and $B_n$ and $n_{test} = n-n_{learn}$ realizations to test the estimator by comparing the
empirical $L_2$-risk of different values of $B_n$ and $K_n$ and choosing the best estimator according to this criterion.  $K_n$ is chosen out of a set 
$ \{4,8,16,32,64,128\}$ and $B_n$ out of a set $\{1,2, \dots, 6, 8, 16,  \dots,  131072\}$.  For each setting of $K_n$ and $B_n$ the estimator is computed ten times with different initializations of the weights and the estimator with the smallest empirical $L_2$ error on the test sample is chosen to compare it with other choices of $K_n$ and $B_n$.
The results of our estimator are compared with standard neural networks,  which are fitted using the adam optimizer in \textit{keras} (\textit{tensorflow} backend) with default learning rate $0.01$ and $1000$ epochs.  In this context we consider structures with one (abbr. \textit{net-1}), three (abbr. \textit{net-3}) and six (abbr. \textit{net-6}) hidden layers.  The number of neurons is also chosen adaptively with the splitting of the sample procedure. As for our estimator we use the set $\{4,8,16,32,64,128\}$ as possible choices for the number of neurons. Furthermore we choose either the ReLU activation function (abbr. \textit{relu-net}) or the sigmoidal activation function (abbr. \textit{sig-net}).
\\
To compare the seven methods (six different network structures with ReLU or sigmoidal activation function + our own method)
we generate $n \in \{200, 400\}$ independent observations from
\begin{align*}
Y= m_i(X) + \sigma_j \cdot \lambda_i \cdot \epsilon \quad (i \in \{1, \dots, 6\},  j \in \{1,2\}),
\end{align*}
where $X$ are uniformly distributed on $[0,1]^d$, $\sigma_j \geq 0$, $\lambda_i \geq 0$ and $\epsilon$ is standard normally distributed and independent of $X$. Thus we use the dataset
\begin{align*}
\mathcal{D}_n = \{(X_1, Y_1), \dots, (X_n,Y_n)\}.
\end{align*}
The value of $\lambda_i$ is chosen in a way that respects the range covered by $m_i$ on the distribution of $X$.  This range is determined empirically as the interquartile range of $10^5$ independent realizations of $m_i(X)$ (and stabilized by taking the median of a hundred repetitions of this procedure), which leads to $\lambda_1=0.24$,  $\lambda_2=0.11$, $\lambda_3=8.76$, $\lambda_4=0.04$, $\lambda_5=0.36$,  and $\lambda_6= 9.11$.
  For the noise value $\sigma_j$ we choose between $5\%$ and $20\%$.
\\
\\
We apply our estimators on the following six regression functions:
\begin{alignat*}{2}
m_1(x) &= \frac{1}{4000} \cdot \sum_{i=1}^7 (x^{(i)})^2-\prod_{i=1}^7\cos\left(\frac{x^{(i)}}{\sqrt{i-1}}\right) \quad & (x \in [0,1]^7)\\
m_2(x) &= \exp\left(0.5 \cdot \sum_{i=1}^{7} (x^{(i)})^2\right), \quad &(x \in [0,1]^{7})\\
m_3(x) &= \sum_{i=1}^{6} 10 \cdot
\left( \left(x^{(i+1)}-(x^{(i)})^2\right)^2+(x^{(i)}-1)^2 \right),
\quad &(x \in [0,1]^7)\\
m_4(x) &= \tanh\left(0.2x^{(1)}+0.9x^{(2)}+x^{(3)} + x{(4)} + 0.2 \cdot x^{(5)} + 0.6 \cdot x^{(6)}\right), \quad &(x \in [0,1]^6)\\
m_5(x) &= \frac{1}{1+\frac{\|x\|}{4}} + (x^{(7)})^2 + x^{(4)} \cdot x^{(5)} \cdot x^{(2)}, \quad &(x \in [0,1]^{10})\\
m_6(x) &= \cot\left(\frac{\pi}{1+\exp\left((x^{(1)})^2+2 \cdot x^{(2)} + \sin(6 \cdot (x^{(4)})^2 -3)\right)}\right) &\\
&\quad + \exp \left(3 \cdot x^{(3)} + 2 \cdot x^{(4)} - 5 \cdot x^{(5)} + \sqrt{x^{(6)} +0.9 \cdot x^{(7)} + 0.1}\right), \quad & (x \in [0,1]^7)
\end{alignat*}

The quality of each of the estimators is determined by the empirical $L_2$-error, i.e.  by 
\begin{align*}
\epsilon_{L_2, N} = \frac{1}{N} \sum_{k=1}^N (m_{n,i}(X_{n+k}) - m_i(X_{n+k}))^2, 
\end{align*}
where $m_{n,i}$ $(i \in \{1,\dots, 6\})$ describes one of the seven estimators based on the $n$ observations and $m_i$ is one of the above mentioned regression functions.  The input values $X_{n+1}, X_{n+2}, \dots, X_{n+N}$ are newly generated independent realizations of the random value $X$. Thus, those values are independent of the $n$ values used for the training and the choice of the parameters of the estimators. We choose $N=10^5$.  Since the value of $\epsilon_{L_2,N}$ strongly depends on the choice of the regression function $m_i$, we normalize this value by dividing it by the 
error of the simplest estimator of $m_i$ , namely the error of a constant function (calculated by the average of the observed data).  The errors in Table \ref{tab1} and \ref{tab2} below are all normalized errors of the form $\epsilon_{L_2,N}(m_{n,i})/\bar{\epsilon}_{L_2,N}(avg)$, where $\bar{\epsilon}_{L_2,N}(avg)$ is the median of $50$ independent realizations one obtains if
one plugs the average of $n$ observations into $\epsilon_{L_2,N}$.  Since our simulation study uses randomly generated data we repeat each estimation $50$ times with different values of $(X, \epsilon)$ in each run.  In the tables below we listed the median (plus interquartile range IQR) of $\epsilon_{L_2,N}(m_{n,i})/\bar{\epsilon}_{L_2,N}(avg)$. 

\begin{table}[h!]
\centering
\begin{tabular}{ccccc}
\hline
\multicolumn{5}{c}{$m_1$}\\
\hline 
\textit{noise} & \multicolumn{2}{c}{$5\%$} & \multicolumn{2}{c}{$20\%$} \\
\hline
\textit{sample size} & $n=200$ & $n=400$ & $n=200$ & $n=400$\\
\hline
$\bar{\epsilon}_{L_2,N}(avg)$ & $0.6823$  & $0.6823$ & $0.6820$  & $0.6825$ \\
\hline
\textit{relu-net-1} & $0.0554 (0.0126)$ & $0.0309 (0.0132)$ & $0.0586 (0.0099)$ & $0.0325 (0.0137)$\\
\textit{relu-net-3} & $0.0632 (0.0110)$ & $0.0428 (0.0100)$ & $0.0645 (0.0094)$ & $0.0449 (0.0135)$\\
\textit{relu-net-6} &$0.0669 (0.0108)$ & $0.0518 (0.0143)$ & $0.0676 (0.0125)$ & $0.0489 (0.0134)$\\
\textit{sig-net-1} & $0.0763 (0.0088)$  & $0.0615 (0.0038)$ & $0.0768 (0.0087)$ & $0.0623 (0.0033)$\\
\textit{sig-net-3} & $0.0967 (0.0150)$  & $0.0683 (0.0036)$ & $0.1006 (0.0140)$ & $0.0685 (0.0042)$\\
\textit{sig-net-6} & $0.1335 (0.0808)$  & $0.0676 (0.0054)$ & $0.1424 (0.0712)$ & $0.0684 (0.0047)$ \\
\textit{comb-classic} & $0.0552 (0.0127)$ & $0.0308 (0.0127)$ & $0.0572 (0.0098)$ & $0.0319 (0.0113)$\\
\textit{lsq-est} & $\mathbf{0.0014} (0.0009)$  & $\mathbf{0.0006} (0.0003)$ & $\mathbf{0.0079} (0.0025)$ & $\mathbf{0.0020} (0.0006)$\\
\textit{comb-new} & $\mathbf{0.0014} (0.0009)$ & $\mathbf{0.0006} (0.0003)$ & $\mathbf{0.0079} (0.0025)$& $\mathbf{0.0020} (0.0006)$ \\
\hline
\end{tabular}
\centering
\begin{tabular}{ccccc}
\hline
\multicolumn{5}{c}{$m_2$}\\
\hline 
\textit{noise} & \multicolumn{2}{c}{$5\%$} & \multicolumn{2}{c}{$20\%$} \\
\hline
\textit{sample size} & $n=200$ & $n=400$ & $n=200$ & $n=400$\\
\hline
$\bar{\epsilon}_{L_2,N}(avg)$ & $2.2190$  & $2.2154$ & $2.2159$ & $2.2151$ \\
\hline
\textit{relu-net-1} & $0.1123 (0.0267)$ & $0.0582 (0.0070)$ & $0.1074 (0.0314)$ & $0.0561 (0.0106)$\\
\textit{relu-net-3} & $0.1015 (0.0252)$ & $0.0607 (0.0104)$ & $0.0977 (0.0266)$ & $0.0629 (0.0144)$\\
\textit{relu-net-6} & $0.0950 (0.0261)$ & $0.0642 (0.0084)$ &$0.0902 (0.0292)$ & $0.0630 (0.0144)$\\
\textit{sig-net-1} & $0.2278 (0.0821)$ & $0.0918 (0.0266)$ &$0.2117 (0.0932)$ & $0.0920 (0.0273)$\\
\textit{sig-net-3} & $0.3058 (0.4057)$  & $0.0745 (0.0117)$ & $0.2487 (0.5568)$ & $0.0736 (0.0123)$\\
\textit{sig-net-6} & $0.9654 (0.6720)$  &$0.1045 (0.0420)$ & $0.8871 (0.6894)$ & $0.0943 (0.0322)$ \\
\textit{comb-classic} & $0.0910 (0.0238)$ & $0.0560 (0.0108)$ & $0.0893 (0.0200)$  & $0.0534 (0.0111)$ \\
\textit{lsq-est} & $\mathbf{0.0303} (0.0095)$  & $\mathbf{0.0167} (0.0043)$ & $\mathbf{0.0294} (0.0113)$ & $\mathbf{0.0169} (0.0052)$\\
\textit{comb-new} & $\mathbf{0.0303} (0.0095)$ & $\mathbf{0.0167} (0.0043)$ & $\mathbf{0.0294} (0.0113)$ & $\mathbf{0.0169} (0.0052)$ \\
\hline
\end{tabular}
\centering
\begin{tabular}{ccccc}
\hline
\multicolumn{5}{c}{$m_3$}\\
\hline 
\textit{noise} & \multicolumn{2}{c}{$5\%$} & \multicolumn{2}{c}{$20\%$} \\
\hline
\textit{sample size} & $n=200$ & $n=400$ & $n=200$ & $n=400$\\
\hline
$\bar{\epsilon}_{L_2,N}(avg)$ & $25.4803$  & $25.4734$  & $25.4906$ & $25.4864$  \\
\hline
\textit{relu-net-1} & $0.6428 (0.2389)$ & $0.2083 (0.0359)$ & $0.7063 (0.1625)$ & $0.2621 (0.0773)$\\
\textit{relu-net-3} & $0.6845 (0.2401)$ & $0.1846 (0.0455)$ & $0.6757 (0.3474)$ & $0.2567 (0.0660)$\\
\textit{relu-net-6} & $0.5204 (0.2810)$ & $0.2114 (0.1316)$ & $0.6455 (0.3179)$ & $0.2878 (0.0960)$\\
\textit{sig-net-1} & $0.9820 (0.0297)$ & $0.9442 (0.0350)$ & $0.9851 (0.0354)$ & $0.9516 (0.0397)$\\
\textit{sig-net-3} & $1.0175 (0.0322)$  & $0.9934 (0.0358)$ & $1.0047 (0.0248)$ & $1.0025 (0.0324)$\\
\textit{sig-net-6} & $1.1014 (0.0275)$  & $1.0016 (0.0253)$ & $1.0064 (0.0271)$ & $1.0058 (0.0247)$ \\
\textit{comb-classic} & $0.4630 (0.2149)$ & $0.1698 (0.0386)$ & $0.5341 (0.2286)$ & $0.2427 (0.0456)$ \\
\textit{lsq-est} & $\mathbf{0.2181} (0.0577)$  & $\mathbf{0.1024} (0.0161)$ & $\mathbf{0.3970} (0.1589)$ & $\mathbf{0.1778} (0.0447)$\\
\textit{comb-new} & $\mathbf{0.2182} (0.0577)$ & $\mathbf{0.1024} (0.0161)$ & $0.4094 (0.1472)$ & $0.1783 (0.0486)$ \\
\hline
\end{tabular}
\caption{Median (and IQR) of the normalized empirical $L_2$-error for each
  of the nine estimates and regression functions $m_1$, $m_2$ and $m_3$ \label{tab1}}
\end{table}
\begin{table}[h!]
\centering
\begin{tabular}{ccccc}
\hline
\multicolumn{5}{c}{$m_4$}\\
\hline 
\textit{noise} & \multicolumn{2}{c}{$5\%$} & \multicolumn{2}{c}{$20\%$} \\
\hline
\textit{sample size} & $n=200$ & $n=400$ & $n=200$ & $n=400$\\
\hline
$\bar{\epsilon}_{L_2,N}(avg)$ & $0.0049$ & $0.0049$ & $0.0049$ & $0.0049$  \\
\hline
\textit{relu-net-1} & $0.2069 (0.1520)$ & $0.0290 (0.0172)$ & $0.1889 (0.1260)$ & $0.0475 (0.0219)$\\
\textit{relu-net-3} & $0.2531 (0.2178)$ & $0.0245 (0.0283)$ & $0.2292 (0.1971)$ & $0.0427 (0.0256)$\\
\textit{relu-net-6} & $0.0760 (0.2178)$ & $0.0276 (0.0290)$ & $0.1058 (0.1280)$ & $0.0349 (0.0386)$\\
\textit{sig-net-1} & $0.1245 (0.0188)$ & $0.0038 (0.0077)$ & $0.1243 (0.0597)$ & $0.0100 (0.0121)$\\
\textit{sig-net-3} & $0.1375 (0.5441)$  & $0.0032 (0.0026)$ & $0.0831 (0.2748)$ & $0.0055 (0.0058)$\\
\textit{sig-net-6} & $0.1766 (0.8157)$  & $0.0147 (0.0284)$  & $0.2949 (0.9867)$ & $0.0117 (0.0195)$ \\
\textit{comb-classic} & $0.0445 (0.0455)$ & $\mathbf{0.0028} (0.0033)$ & $0.0594 (0.0620)$ & $\mathbf{0.0052} (0.0067)$\\
\textit{lsq-est} & $\mathbf{0.0188} (0.0106)$  & $0.0077 (0.0046)$ & $0.0469 (0.0187)$  & $0.0230 (0.0079)$\\
\textit{comb-new} & $\mathbf{0.0184} (0.0142)$ & $\mathbf{0.0028} (0.0033)$ & $\mathbf{0.0437} (0.0262)$ & $\mathbf{0.0052} (0.0050)$ \\
\hline
\end{tabular}
\centering
\begin{tabular}{ccccc}
\hline
\multicolumn{5}{c}{$m_5$}\\
\hline 
\textit{noise} & \multicolumn{2}{c}{$5\%$} & \multicolumn{2}{c}{$20\%$} \\
\hline
\textit{sample size} & $n=200$ & $n=400$ & $n=200$ & $n=400$\\
\hline
$\bar{\epsilon}_{L_2,N}(avg)$ & $0.0135$  & $0.0135$  & $0.1345$   &  $0.0135$ \\
\hline
\textit{relu-net-1} & $0.1989 (0.0737)$ & $0.1019 (0.0264)$ & $0.2846 (0.1094)$ & $0.1768 (0.0618)$\\
\textit{relu-net-3} & $0.1826 (0.0667)$ & $0.1061 (0.0454)$ & $0.3182 (0.1858)$ & $0.2230 (0.0808)$\\
\textit{relu-net-6} & $0.1744 (0.1099)$ & $0.1066 (0.0264)$ & $0.3519 (0.1995)$ & $0.2005 (0.0621)$\\
\textit{sig-net-1} & $0.0971 (0.0183)$ & $0.0784 (0.0071)$ & $0.1139 (0.0202)$ & $0.0880 (0.0131)$\\
\textit{sig-net-3} & $0.0935 (0.0202)$  & $0.0796 (0.0049)$ & $0.1190 (0.0434)$ & $0.0880 (0.0138)$\\
\textit{sig-net-6} & $0.1396 (0.0547)$  & $0.0850 (0.0087)$ & $0.1816 (0.0434)$ & $0.0983 (0.0207)$\\
\textit{comb-class} & $0.0948 (0.0187)$ & $0.0794 (0.0869)$ & $\mathbf{0.1121} (0.0210)$ & $\mathbf{0.0870} (0.0116)$\\
\textit{lsq-est} & $0.1025 (0.0343)$ & $\mathbf{0.0516} (0.0118)$  & $0.1682 (0.0434)$ & $0.0981 (0.0181)$\\
\textit{comb-new} & $\mathbf{0.0923} (0.0198)$ & $\mathbf{0.0516} (0.0118)$ & $0.1162 (0.0392)$& $0.0918 (0.0131)$ \\
\hline
\end{tabular}
\centering
\begin{tabular}{ccccc}
\hline
\multicolumn{5}{c}{$m_6$}\\
\hline 
\textit{noise} & \multicolumn{2}{c}{$5\%$} & \multicolumn{2}{c}{$20\%$} \\
\hline
\textit{sample size} & $n=200$ & $n=400$ & $n=200$ & $n=400$\\
\hline
$\bar{\epsilon}_{L_2,N}(avg)$ & $591.77$   & $592.74$& $591.10$ & $589.81$  \\
\hline
\textit{relu-net-1} & $0.1897 (0.0913)$ & $0.0552 (0.0346)$  & $0.2195 (0.0876)$  & $0.0845 (0.0378)$ \\
\textit{relu-net-3} & $0.1349 (0.0828)$ & $0.0398 (0.0315)$  & $0.1383 (0.0828)$  & $0.0467 (0.0363)$  \\
\textit{relu-net-6} & $0.1075 (0.0856)$  & $0.0385 (0.0247)$  & $\mathbf{0.1073} (0.0657)$ & $0.0434 (0.0249)$  \\
\textit{sig-net-1} & $0.5733 (0.1122)$ & $0.2027 (0.0539)$ & $0.5572 (0.0852)$ & $0.1877 (0.1123)$\\
\textit{sig-net-3} & $0.8355 (0.1451)$& $0.3905 (0.1181)$ & $0.8697 (0.1728)$ & $0.3805 (0.0822)$\\
\textit{sig-net-6}  &$1.0058 (0.0728)$ & $0.7017 (0.1435)$ & $0.9955 (0.0912)$ & $0.7332 (0.1188)$\\
\textit{comb-classic} & $\mathbf{0.0991} (0.0749)$ & $\mathbf{0.0331} (0.0242)$ & $0.1173 (0.0842)$ & $\mathbf{0.0352} (0.0229)$\\
\textit{lsq-est} & $0.3070 (0.1534)$ & $0.1784 (0.0758)$ & $0.2643 (0.1951)$ & $0.1655 (0.0570)$\\
\textit{comb-new} & $0.1063 (0.0803)$ & $\mathbf{0.0331} (0.0394)$ & $0.1214 (0.0902)$ & $\mathbf{0.0352} (0.0229)$\\
\hline
\end{tabular}
\centering
\caption{Median (and IQR) of the normalized empirical $L_2$-error for each
  of the nine estimators and regression functions
  $m_4$, $m_5$ and $m_6$ \label{tab2}}
\end{table}
We observe that our new linear least squares estimator outperforms the other approaches in 15 of 24 cases.  Especially, for the functions $m_1, m_2$ and $m_3$ our estimator is always the best and has, as for function $m_1$,  a more than 10 times smaller error than the error of the second best approach.  For this function we also observe that the relative improvement of our estimator with an increasing sample size is often much larger than the improvement of the other approaches.  This can be considered as an indicator for a better rate of convergence of the estimator.
\\
With regard to the other three functions $m_4, m_5$ and $m_6$ our estimator is only sometimes the best.  Especially for the cases of $m_4$ and $m_6$ the results are not entirely satisfactory.  With regard to $m_6$ our estimator is always outperformed by the standard ReLU networks with six hidden layers and also for $m_4$ and the higher sample size the standard sigmoidal networks with one or three hidden layers are better by a factor of at least two.
\\
With the goal in mind to construct an estimator based on statistical theory that provides satisfactory results in all settings, we have extended our simulation study.  We constructed a combined estimator (abbr. \textit{comb-new}) , i.e. , an estimator which chooses between the new least squares estimator and the standard nets the one with the smallest empirical $L_2$-error on the dataset $x_{test}$. This estimator was compared to a classical combined estimator (abbr.  \textit{comb-classic}),  i.e.  an estimator that chooses the best standard net according to the smallest empirical $L_2$ error on the test sample. The results are also given in Table \ref{tab1} and \ref{tab2}.
In $16$ of $24$ cases our new combined estimator is better than the classical approach. In four cases both estimators are of the same size and only in the remaining four cases the classical approach is slightly better.  For function $m_1$ our new combined estimator is more than $15-20$ times better than the classical approach and also in most of the other cases we see a significant difference between the error of the new combined estimator and the classical one.  A look at the results in which our estimator performs somewhat worse shows that the classical estimator is always less than $10\%$ better.  For us these small changes of the median error are not significant.

Summarizing our simulation study we see that
our newly proposed combined estimator
is in our simulation study
never significantly worse than the standard neural network estimators, but
is in
some of the considered cases much better than the standard estimators.
\\
\\
Since the values of $K_n$ and $B_n$ that define our new least squares estimator are chosen by a splitting of the sampling procedure, it is of particular interest how sensitive the estimator is to different choices of $K_n$ and $B_n$.  For a slightly reduced set of possible $K_n$ and $B_n$, i.e. , a set of $[4, 16, 64]$ for $K_n$ and $[1,16,64,256, 1024]$ for $B_n$, the next plots show how different the empirical $L_2$ error values behave for the estimator.
%

\begin{figure}[h!]
 \centering
 \subfigure[]{\includegraphics[width=0.3\textwidth]{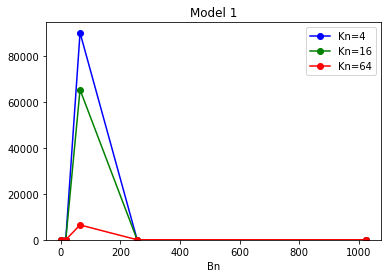}}
  \subfigure[]{\includegraphics[width=0.3\textwidth]{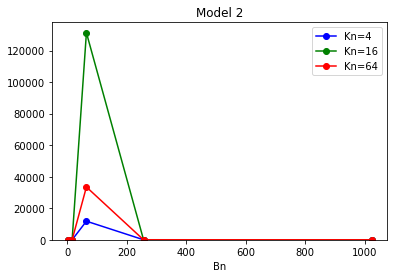}}
    \subfigure[]{\includegraphics[width=0.3\textwidth]{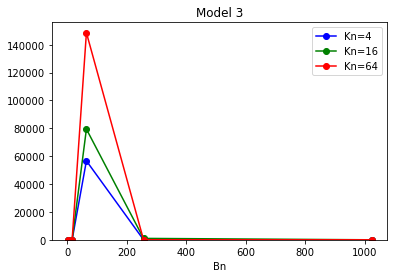}}
     \subfigure[]{\includegraphics[width=0.3\textwidth]{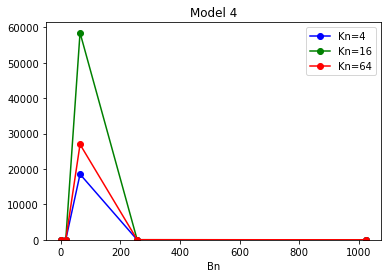}}
  \subfigure[]{\includegraphics[width=0.3\textwidth]{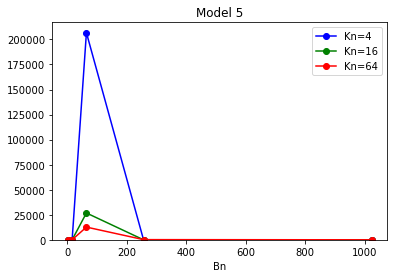}}
    \subfigure[]{\includegraphics[width=0.3\textwidth]{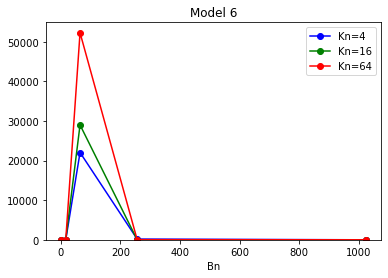}}
  \caption{Comparison of different settings of $B_n$ and $K_n$ for all models, noise $5\%$ and sample size $n=200$}
  \label{fig:bn}
\end{figure}

For all possibilities of $K_n$, the estimator has a very high error for $B_n = 64$.  For small $B_n$ of $1$ or $16$ and for high $B_n$ of $256$ or $1024$ the error is closer to zero. A closer look on the numbers shows that the error nevertheless varies between $0$ and $15$ with a tendency for smaller values of $B_n$ to show better performance.  For different models, different values of $K_n$ are the best choice.  This shows that the splitting of the sample procedure is indeed important for the performance of the estimator.

\section{Proofs}
\label{se5}

\subsection{Approximation error of shallow neural networks
  with random inner weights}

The main trick in our results is the following approximation
result for shallow neural networks with random inner weights
and threshold squasher $\sigma(x)=1_{[0,\infty)}(x)$
  as activation function.
 The proof relies on an extension of the approach in
Barron (1994).  Here condition (\ref{le11eq5}) below will be used to show that for the
  logistic squasher $\sigma(x)=1/(1+\exp(-x))$ we have
  \[
  \sigma \left(
B_n \cdot (W_k^T X_i + T_k)
\right)
\approx
1_{[0,\infty)} \left(
W_k^T X_i + T_k
\right).
  \]
 
\begin{lemma}
  \label{le11}
  Let $r>0$ and $\tilde{K}_n \in \N$ with $\tilde{K}_n \leq n$,
  and set
  $K_n = (\lceil \log n \rceil)^4 \cdot \tilde{K}_n$.
  Let $m:\Rd \rightarrow \R$ be a function with
  \begin{equation}
    \label{le11eq1}
\int_{\Rd} | m(x)| \, dx < \infty,
  \end{equation}
  and assume that
  the Fourier transform $\hat{F}$ of $m$ satisfies
  \begin{equation}
    \label{le11eq2}
        |\hat{F}(\omega)| \leq
    \frac{c_1}{\| \omega \|^{d+1} \cdot (\log \| \omega\|)^2 }
    \quad
    \mbox{for all } \omega \in \Rd \mbox{ with }
    \| \omega \| \geq 2
  \end{equation}
  for some $c_1>0$.
  Let $(X,Y)$, $(X_1,Y_1)$, \dots, $(X_n,Y_n)$ be $\Rd \times \R$--valued
  random variables, and let $W_1$, \dots, $W_{K_n}$, $T_1$,
  \dots, $T_{K_n}$ be independent random variables, independent
  from  $(X,Y)$, $(X_1,Y_1)$, \dots, $(X_n,Y_n)$, such that
  $W_1$, \dots, $W_{K_n}$ are uniformly distributed
  on $\{ x \in \Rd \, : \, \|x\|=1 \}$ and \linebreak
   $T_1$,
  \dots, $T_{K_n}$ are uniformly distributed on $[-r,r]$.
  Then
for $n$ sufficiently large,
  there exist (random)
  \[
  \alpha_0 \in [-c_3, c_3] \quad \mbox{and} \quad
  \alpha_1, \dots, \alpha_{K_n} \in \left[
- \frac{c_4}{\tilde{K}_n}, \frac{c_4}{\tilde{K}_n}
    \right],
  \]
  which are independent of $(X,Y)$, $(X_1,Y_1)$, \dots, $(X_n,Y_n)$,
  such that outside of an event with probability less than
  or equal to
  \begin{equation}
    \label{le11eq3}
  \exp \left(-\frac{1}{4} \cdot (\lceil \log n \rceil)^2 \right)
  +
  \frac{1}{n},
  \end{equation}
  we have
  \begin{equation}
    \label{le11eq4}
  \int_{S_r}
  \left|
  m(x)-
  \alpha_0
  -
  \sum_{k=1}^{K_n} \alpha_k \cdot
  1_{[0,\infty)} \left(
W_k^T \cdot x + T_k
    \right)
    \right|^2
    \PROB_X (dx)
    \leq
    \frac{c_{22}}{\tilde{K}_n}
  \end{equation}
  and
  \begin{equation}
    \label{le11eq5}
      \min_{i=1, \dots, n, k=1, \dots, K_n: \atop \alpha_k \neq 0}
      |         W_k^T X_i +T_k| \geq \delta_n,
      \end{equation}
   where $\delta_n=\frac{r}{n^{2} \cdot K_n}$.
  \end{lemma}

\begin{proof}
The complete proof of this result is found in the supplement.
\end{proof}

\subsection{Outline of the proof of Theorem \ref{th1}}
In this subsection we give an outline of the proof of
Theorem \ref{th1}. The complete proof is given in the
supplement.

W.l.o.g. we assume $\|m\|_\infty \leq \kappa_n$.
Set $\tilde{K}_n = \lfloor K_n / (\lceil \log n \rceil)^4 \rfloor$ and
let $A_n$ be the event that
$|Y_i| \leq \kappa_n$ holds for all $i \in \{1, \dots, n\}$ and that
there exist
(random)
  \begin{align}
  \label{alpha}
  \alpha_0 \in [-c_3, c_3] \quad \mbox{and} \quad
  \alpha_1, \dots, \alpha_{K_n} \in \left[
- \frac{c_4}{\tilde{K}_n}, \frac{c_4}{\tilde{K}_n}
    \right],
  \end{align}
  which are independent of $(X,Y)$, $(X_1,Y_1)$, \dots, $(X_n,Y_n)$,
  such that (\ref{le11eq4}) and
  \begin{equation}
    \label{pth1outeq**}
    \min_{i=1, \dots, n, k=1, \dots, K_n: \atop
      \alpha_k \neq 0}
    \left|
(\beta_k^{(0)})^T X_i + \gamma_k^{(0)}
    \right|
    \geq \delta_n
    \end{equation}
  hold for
  $\delta_n=B_n \cdot \frac{\sqrt{d}}{n^{2} \cdot K_n}$,
  i.e., for $\delta_n=(\log n)^2$.

We have
\begin{eqnarray*}
  &&
  \EXP \int | m_n(x)-m(x)|^2 \PROB_X (dx)
  \\
  &&
  =
  \EXP \left( \int | m_n(x)-m(x)|^2 \PROB_X (dx) \cdot 1_{A_n} \right)
  +
    \EXP \left( \int | m_n(x)-m(x)|^2 \PROB_X (dx) \cdot 1_{A_n^c} \right)
    \\
    &&
    \leq
      \EXP \left( \int | m_n(x)-m(x)|^2 \PROB_X (dx) \cdot 1_{A_n} \right)
  +
4 \kappa_n^2 \cdot \PROB(A_n^c)
\\
&&
=
\mathbf E \Bigg(\bigg( \int |m_n(x) - m(x)|^2 \PROB_X (dx)
  \\
  &&
  \quad \quad \quad
  -
  2 \cdot \bigg( \frac{1}{n} \sum_{i=1}^n |\tilde{m}_n(X_i)-Y_i|^2
  \cdot 1_{\{|Y_j| \leq \kappa_n \, (j \in \{1, \dots, n\} ) \}}
  -
  \frac{1}{n} \sum_{i=1}^n |m(X_i)-Y_i|^2
  \bigg)
  \bigg)
  \cdot 1_{A_n}
  \Bigg)
  \\
  &&
  \quad
  +
  2 \cdot
  \EXP \left(
\bigg( \frac{1}{n} \sum_{i=1}^n |\tilde{m}_n(X_i)-Y_i|^2
  \cdot 1_{\{|Y_j| \leq\kappa_n \, (j \in \{1, \dots, n\} ) \}}
  -
  \frac{1}{n} \sum_{i=1}^n |m(X_i)-Y_i|^2
  \bigg)
  \cdot 1_{A_n}
  \right)
  \\
  &&
  \quad
  +
  4 \kappa_n^2 \cdot \PROB(A_n^c)
  \\
  &&
  =: T_{1,n}+T_{2,n}+T_{3,n}.
\end{eqnarray*}
Using results from empirical process theory together
with bounds on the norm of the weights occuring
during gradient descent we  show in the supplement
\[
  T_{1,n} \leq c_{29} \cdot \frac{(\log n)^3 \cdot K_n}{n}.
\]
Furthermore we will use Lemma \ref{le11} to show
\[
T_{3,n}
\leq
c_{30} \cdot
\frac{(\log n)^2}{n}.
\]
The remaining term we have to bound is $T_{2,n}$.
Let
$\alpha_0$ \dots $\alpha_{K_n}$ be defined as in \eqref{alpha}
and
define on $[0,1]$ a piecewise constant approximation of $m$ by
\[
f(x) =
\alpha_0 +
  \sum_{k=1}^{K_n} \alpha_k \cdot
  1_{[0,\infty)} \left(
\sum_{j=1}^d w_{k,j}^{(0)} \cdot x^{(j)} + w_{k,0}^{(0)}
    \right).
  \]
Set 
\[ 
	f^*(x) 
	=
        \alpha_0 +
        \sum_{k=1}^{K_n} \alpha_k \cdot
  \sigma \left(
\sum_{j=1}^d w_{k,j}^{(0)} \cdot x^{(j)} + w_{k,1}^{(0)}
    \right).
\]
For $g(x)= \alpha_0 + \sum_{k=1}^{K_n} \alpha_k \cdot \sigma(\beta_k^T \cdot x + \gamma_k)$
we define
\[
\text{pen}(g) 
= 
\frac{c_2}{K_n} \cdot \sum_{k=0}^{K_n} (\alpha_k)^2.
\]
We will see
that it is enough to show
\begin{eqnarray*}
  &&
	F( \alpha^{(t_n)} , \beta^{(t_n)} , \gamma^{(t_n)})
	-
	\frac{1}{n} \sum_{i=1}^n |Y_i - f^{*}(X_i)|^2 -\text{pen}(f^*)
        \\
        &&
        \leq
c_{31} \cdot \frac{(\log n)^3 \cdot K_n}{n},  
\end{eqnarray*}
where
\[
F(\alpha,\beta,\gamma)
=
\frac{1}{n} \sum_{i=1}^n |Y_i - f_{net,(\alpha,\beta,\gamma)}(X_i)|^2
+
\frac{c_2}{K_n}
 \cdot \sum_{k=0}^{K_n} (\alpha_k)^2.
\]
The main trick is to deduce via an elementary gradient descent 
analysis that on $A_n$ 
\begin{eqnarray*}
&&
	F( \alpha^{(t+1)} , \beta^{(t+1)} , \gamma^{(t+1)})
	-
	\frac{1}{n} \sum_{i=1}^n |Y_i - f^{*}(X_i)|^2 -{\text pen}(f^*)
         \\
&& \leq
\left(
1 - \frac{2 \cdot c_2}{L_n \cdot K_n}
\right)
\cdot
\left(
	F( \alpha^{(t)} , \beta^{(t)} , \gamma^{(t)})
	-
	\frac{1}{n} \sum_{i=1}^n |Y_i - f^{*}(X_i)|^2 -{\text pen}(f^*)
\right)
        \\
&&
\quad
+
\frac{2 \cdot c_2}{L_n \cdot K_n} \cdot
\left(
F( \alpha^{*} , \beta^{(0)} , \gamma^{(0)})
-
F( \alpha^{*} , \beta^{(t)} , \gamma^{(t)})
\right),
\end{eqnarray*}
for any $t \in \{1, \dots, t_n-1\}$ and $\alpha^*=(\alpha_k)_{k=0, \dots, K_n}$ with $\alpha_{k}$
as defined above.  If we now take advantage of the fact that in the logistic squasher the inner weights (which are not equal to zero) 
change only slightly, this implies
\begin{eqnarray*}
&&
	F( \alpha^{(t_n)} , \beta^{(t_n)} , \gamma^{(t_n)})
	-
	\frac{1}{n} \sum_{i=1}^n |Y_i - f^{*}(X_i)|^2 -\text{pen}(f^*)
        \\
        &&
        \leq
\left(
1 - \frac{2 \cdot c_2}{L_n \cdot K_n}
\right)^{t_n}
\cdot
\left(
	F( \alpha^{(0)} , \beta^{(0)} , \gamma^{(0)})
	-
	\frac{1}{n} \sum_{i=1}^n |Y_i - f^{*}(X_i)|^2 -\text{pen}(f^*)
\right)
   \\
   &&
   \quad
   +
   t_n \cdot
   \frac{2 \cdot c_2}{L_n \cdot K_n}
   \cdot
 c_{32} \cdot (\kappa_n + (\log n)^4) \cdot
 c_{33} \cdot
 \frac{1}{n}
\\
 &&
 \leq \exp \left(
- \frac{2 \cdot c_2}{L_n \cdot K_n} \cdot t_n
\right)
\cdot
F( \alpha^{(0)} , \beta^{(0)} , \gamma^{(0)})
+ c_{34} \cdot \frac{\kappa_n + (\log n)^4}{n}
\\
&&
\leq
\exp( - c_{35} \cdot (\log n)^2 ) \cdot c_{36} \cdot (\log n)^2
+ c_{37} \cdot \frac{\kappa_n+(\log n)^4}{n}
\leq
 c_{38} \cdot \frac{(\log n)^4}{n}.
  \end{eqnarray*}

\hfill $\Box$

\subsection{Proof of Corollary \ref{co1}}
It is well known that Lebesgue integrability and radial symmetry of $m$ imply radial symmetry of $|\hat{F}|$
(cf., Theorem 4.5.3 in Epstein (2008)). Then in Lemma \ref{le11} one can replace
(\ref{le11eq2})
by \eqref{c1eq1}, where in the proof the radial symmetric density $g$ is defined by 
        \begin{align*}
        g(t, \omega) = \begin{cases}
        \frac{1}{2r} \cdot \mathds{1}_{[-r,r]}(t) \cdot c_{39}, \quad &\text{if} \ \|\omega\| \leq 1\\
        \frac{1}{2r} \cdot \mathds{1}_{[-r,r]}(t) \cdot c_{39} \cdot \|\omega\| \cdot |\hat{F}(\omega)|, \quad &\text{if} \ \|\omega\| > 1
        \end{cases}
        \end{align*}
        with a suitable constant $c_{39}>0$. Now the proof of Corollary 1 is analogous to the proof of Theorem \ref{th1}. 
\hfill $\Box$

\subsection{Outline of the proof of Theorem \ref{th3}}
We set $p=d/2+1$, $M_n= \lceil n^{1/(2p+d)} \rceil$ and define
\begin{align*}
  m^{(c_n)} (x) = \sum_{j=1}^{M_n^d} c_{n,j} \cdot (\log n)^{-3} \cdot
  M_n^{-p} g(M_n(x-a_{n,j}))
\end{align*}
as in the proof of Theorem 3.2 in Gy\"orfi et al. (2002).
Here we choose $g$ so smooth that
\begin{align*}
  |\hat{F}_g(\omega)| \leq
      \frac{c_{40}}{\| \omega \|^{d+1} \cdot (\log \| \omega\|)^2 }
    \quad
    \mbox{for all } \omega \in \Rd \mbox{ with }
    \| \omega \| \geq 2
\end{align*}
holds (cf., Remark \ref{re3}). Let $\mathcal{C}_n$ be the set of all $(c_{n,j})_{j \in \{1, \dots, M_n^d\}} \in \{-1,1\}^{M_n^d}$, where 
\begin{align*}
\left|\sum_{j=1}^{M_n^d} c_{n,j} \cdot \exp(i \omega	^T a_{n,j}) \right| \leq (\log n) \cdot M_n^{d/2} \quad (\omega \in \Rd \setminus \{0\})
\end{align*}
holds.  If $c_n = (c_{n,j})_j \in \mathcal{C}$ we can show that
\begin{align*}
  |\hat{F}_{m^{(c_n)}}(\omega)|  \leq
    \frac{c_1}{\| \omega \|^{d+1} \cdot (\log \| \omega\|)^2 }
    \quad
    \mbox{for all } \omega \in \Rd \mbox{ with }
    \| \omega \| \geq 2
.
\end{align*}
Let $X, N, C_n$ be independent with $X\sim U([0,1]^d)$, $N \sim \mathcal{N}(0,1)$ and $C_n \sim U(\mathcal{C}_n)$.
Arguing as in the proof of
Theorem 3.2 in Gy\"orfi et al. (2002) we can bound 
\begin{align*}
&\inf_{\hat{m}_n} \sup_{(X,Y) \in \mathcal{D}} \EXP \int |\hat{m}_n(x) - m(x)|^2 \PROB_X(dx)\\
&\geq c_{41} \cdot (\log n)^{-6} \cdot n^{-\frac{2p}{2p+d}} - 4 \cdot \|g\|_\infty^2 \cdot \PROB\{C_n \notin \mathcal{C}_n\}.
\end{align*}
The result follows from the definition of $p$ and an application of
Hoeffding's inequality, which yields
\begin{align*}
\PROB\{C_n \notin \mathcal{C}_n\} \leq \frac{c_{42}}{n}.
\end{align*}
\quad \hfill $\Box$

\subsection{Outline of the proof of Theorem \ref{th2}}
  
W.l.o.g. we assume $\|m\|_\infty \leq \kappa_n $.
Set $\tilde{K}_n = \lfloor K_n / (\lceil \log n \rceil)^4 \rfloor$ and
let $A_n$ be the event that
$|Y_i| \leq \kappa_n$ holds for all $i=1, \dots, n$ and that
there exist
(random)
  \[
  \alpha_0 \in [-c_3, c_3] \quad \mbox{and} \quad
  \alpha_1, \dots, \alpha_{K_n} \in \left[
- \frac{c_4}{\tilde{K}_n}, \frac{c_4}{\tilde{K}_n}
    \right],
  \]
  which are independent of $(X,Y)$, $(X_1,Y_1)$, \dots, $(X_n,Y_n)$,
  such that (\ref{le11eq4}) and
  (\ref{pth1outeq**}) hold for
  $\delta_n=B_n \cdot \frac{\sqrt{d}}{n^{2} \cdot K_n}$,
  i.e., for $\delta_n=(\log n)^2$. Furthermore, define
  $f^*$ as in the proof of Theorem \ref{th1}.

As in the proof of Theorem \ref{th1} we get
\begin{eqnarray*}
  &&
  \EXP \int | m_n(x)-m(x)|^2 \PROB_X (dx)
\\
&&
\leq
\mathbf E \Bigg(\bigg( \int |m_n(x) - m(x)|^2 \PROB_X (dx)
  \\
  &&
  \quad
  -
  2 \cdot \bigg( \frac{1}{n} \sum_{i=1}^n |\tilde{m}_n(X_i)-Y_i|^2
  \cdot 1_{\{|Y_j| \leq \kappa_n \, (j \in \{1, \dots, n\} ) \}}
  -
  \frac{1}{n} \sum_{i=1}^n |m(X_i)-Y_i|^2
  \bigg)
  \bigg)
  \cdot 1_{A_n}
  \Bigg)
  \\
  &&
  \quad
  +
  2 \cdot
  \EXP \left(
\bigg( \frac{1}{n} \sum_{i=1}^n |\tilde{m}_n(X_i)-Y_i|^2
  \cdot 1_{\{|Y_j| \leq \kappa_n \, (j \in \{1, \dots, n\} ) \}}
  -
  \frac{1}{n} \sum_{i=1}^n |m(X_i)-Y_i|^2
  \bigg)
  \cdot 1_{A_n}
  \right)
  \\
  &&
  \quad
  +
  4 \kappa_n^2 \cdot \PROB(A_n^c)
  \\
  &&
  =: T_{1,n}+T_{2,n}+T_{3,n}.
\end{eqnarray*}
Standard results from empirical process theory enable us to show
\[
T_{1,n}
\leq
 \frac{c_{43} \cdot
                                                  (\log n)^2\cdot
   \left(
   (K_n+1) \cdot (\log n)
 +1
\right)
 }{n}
 \leq
 c_{44} \cdot 
 \frac{(\log n)^3 \cdot K_n}{n}.
\]
As in the proof of Theorem \ref{th1} application
of Lemma \ref{le11} yields
\[
T_{3,n}
\leq
c_{45} \cdot
\frac{(\log n)^2}{n}.
\]
Hence again it will be sufficient to derive a bound
on $T_{2,n}$, and as in the proof of Theorem \ref{th1}
the crucial step will be to bound on the event $A_n$
\[
\frac{1}{n} \sum_{i=1}^n |\tilde{m}_n(X_i)-Y_i|^2
  -
  \frac{1}{n} \sum_{i=1}^n |f^*(X_i)-Y_i|^2.
\]
But this is quite simple here as
by definition of $\tilde{m}_n$ as a least
squares estimator this term is less
than or equal to zero.
\hfill $\Box$


\section*{Acknowledgements}
The authors are grateful to the two anonymous referees and
the Associate Editor Mark Podolskij for their constructive comments that improved the
quality of this paper.



\begin{supplement}
\stitle{Appendix A: Further proofs}
\sdescription{Appendix A contains the complete proofs of
  Theorem \ref{th1}, Theorem \ref{th3}, Theorem \ref{th2}
and Corollary \ref{cor31}.}
\end{supplement}
\begin{supplement}
\stitle{Appendix B: Further simulation results}
\sdescription{Appendix B contains the simulation results of Section \ref{se4} for sample size $N=1000$.}
\end{supplement}


\begin{thebibliography}{4}



\bibitem{AlLiLi20}
  Allen-Zhu, Z.,  Li, Y.,  and Liang, Y. (2019).  Learning and
  generalization in overparameterized neural networks, going beyond two layers.  In
{\it Advances in Neural Information Processing Systems}, {\bf 32}.

    
  \bibitem{AlLiSo19}
  Allen-Zhu, Z., Li, Y., and Song, Z. (2019).
  A convergence theory for deep learning via over-parameterization. In 
  \textit{International Conference on Machine Learning},
  {\bf 97}, pp. 242-252.

  \bibitem{APCZ}
  Andoni, A., Panigraphy, R., Valiant, G., and Zhang, L.(2014).
  Learning polynomials with neural networks. In
  \textit{International Conference on Machine Learning}, pp. 1908 --1916.

 \bibitem{ArCoGoHu18}
 Arora, S., Cohen, N., Golowich, N., and Hu, W. (2018).
 A convergence analysis of gradient descent for deep linear neural networks. In
 \textit{International Conference on Learning Representations}.


  \bibitem{ArDuHuLiSaWa19}
 Arora, S., Du, S., Hu, W., Li, Z., Salakhutdinov, R., and Wang, R. (2019a).
 On exact computation with an infinitely wide neural net. In
 \textit{Advances in Neural Information Processing Systems}, {\bf 32}.
 
 
\bibitem{ArDuHuLiWa}
Arora, S., Du, S., Hu, W., Li, Z., and Wang, R. (2019b).
Fine-grained analysis of optimization and generalization for overparameterized two-layer neural networks. In 
\textit{International Conference on Machine Learning},
{\bf 97}, pp. 322-332.


  \bibitem{BaClKo09}
Bagirov, A. M., Clausen, C., and Kohler, M. (2009).
Estimation of a regression function by maxima of minima of linear functions.
{\it IEEE Transactions on Information Theory}
{\bf 55}: 833-845.


\bibitem{Bar93}
Barron, A. R. (1993).
Universal approximation bounds for superpositions of
a sigmoidal function.
{\it IEEE Transactions on Information Theory} {\bf 39}: 930--944.

  
\bibitem{Bar94}
Barron, A. R. (1994).
Approximation and estimation bounds for artificial
neural networks.
{\it Machine Learning} {\bf 14}: 115-133.
   

\bibitem{BHLM}
   Bartlett, P., Harvey, N., Liaw, C., and Mehrabian, A. (2019).
   Nearly-tight VC-dimension bounds for piecewise linear neural networks.
   {\it Journal of Machine Learning Research} {\bf 20}:1--17.
   


\bibitem{BK19}
  Bauer, B., and Kohler, M. (2019).
On deep learning as a remedy for the curse of dimensionality in nonparametric regression. {\it Annals of Statistics} {\bf 47}: 2261--2285.

\bibitem{BrKoWa19}
  Braun, A., Kohler, M., and Walk, H. (2019).
  On the rate of convergence of a neural network regression estimate learned by gradient descent. Preprint, arXiv: 1912.03921.
  

\bibitem{BGMS}
Brutzkus, A., Globerson, A., Malach, E., and Shalew-Shwartz S. (2018).
SGD learn overparametrized networks that provably generalize on linearly seperable data. In
\textit{International Conference on Learning Representation}



\bibitem{CHMAL15}
  Choromanska, A., Henaff, M., Mathieu, M., Arous, G. B., and LeCun, Y.
  (2015)
  The loss surface of multilayer networks.  In
  \textit{International  Conference
  on Artificial Intelligence and Statistics}, {\bf 38}, pp. 192-204.




\bibitem{DuLe18}
Du, S., and Lee, J. (2018).
On the power of over-parametrization in neural networks with quadratic activation. In
\textit{International Conference on Machine Learning}, 
{\bf 80}, pp. 1329-1338.


\bibitem{Duetal18}
  Du, S., Lee, J., Li, H., Wang, L., und Zhai, X. (2019).
  Gradient descent finds global minima of deep neural networks. In
  \textit{International Conference on Machine Learning}, pp. 1675-1685.
 
 
\bibitem{Dudek}
Dudek, G. (2019). Generating random weights and biases in feedfoward neural networks with random hidden nodes. 
\textit{Information Sciences}, {\bf 481}(C): 33-56.


\bibitem{Epst03}
  Epstein, Ch.~L. (2008).
  {\it Introduction to the Mathematics of Medical Imaging}.
  2nd edition, SIAM. Society for Industrial and Applied Mathematics.
  Philadelphia.
  \bibitem{GMMM19}
  Ghorbani, B., Mei, S., Misiakiewicz, T. and Montanari, A. (2019). Limitations of lazy training of two-layer neural networks. In 
  \textit{Advances in Neural Information Processing Systems}, {\bf 32}.
  
  \bibitem{Gonon}
  Gonon, L. (2021).
  Random feature neural networks learn Black-Scholes type PDEs without curse of dimensionality. Preprint, {\it arXiv: 2106.08900}.
  
  \bibitem{Gonon2}
  Gonon, L., Grigoryeva, L., and Ortega, J.-P. (2023). Approximation bounds for random neural networks and reservoir systems. 
  \textit{The Annals of Applied Probability}, {\bf 33}(1): 28--69.
  
 
\bibitem{GoBeCo2016}
  Goodfellow, I., Bengio, Y., and Courville, A.
  (2016).
  {\it Deep Learning}.
  MIT Press, Cambridge, Massachusetts.
  
 \bibitem{GKKW02}
 Gy\"orfi, L., Kohler, M., Krzy\.zak, A., and Walk, H. (2002).
 {\it A Distribution--Free Theory of Nonparametric Regression}.
 Springer, New York.
 
 \bibitem{Huang}
 Huang, G.-B., Chen, L., and Siew, C.-K. (2006). Universal approximation using incremental contractive feedforward networks with random hidden nodes. {\it IEEE Transactions on Neural Networks} {\bf 17}: 879-892.
 
 \bibitem{IgPa95}
 Igelnik, B. and Pao, Y.-H. (1995).
 Stochastic choice of basis functions in adaptive function approximation and the functional-link net.
 {\it IEEE Transactions on Neural Networks}, {\bf 6}: 1320 -1329.
 
 
 
\bibitem{ImFu18}
Imaizumi, M., and Fukamizu, K. (2018).
Deep neural networks learn non-smooth functions effectively. In
{\it International Conference on Artificial Intelligence and Statistics}, pp. 869-878.


\bibitem{JaGaHo20}
  Jacot, A., Gabriel, F., und Hongler, C. (2020).
  Neural Tangent Kernel: Convergence and Generalization in
  Neural Networks.  In 
  {Advances in Neural Information Processing Systems}, {\bf 31}.

\bibitem{JaMoMo21}
  Javanmard, A., Mondelli, M., and Montanari, A. (2021).
  Analysis of two-layer neural network via displacement
  convexity.
 {\it Annals of Statistics}, {\bf 48}(6): 3619-3642.

\bibitem{Ka16}
Kawaguchi, K. (2016).
Deep learning without poor local minima. 
\textit{Advances in Neural Information Processing Systems}, {\bf 29}.



\bibitem{KaHu19}
Kawaguchi, K, and Huang, J. (2019).
Gradient descent finds global minima for  generalizable deep neural networks of practical sizes. In
\textit{2019 57th annual allerton conference on communication, control, and computing (Allerton)}, pp. 92--99.



\bibitem{Kim14}
Kim, Y. (2014).
Convolutional Neural Networks for Sentence Classification. In 
\textit{Empirical Methods in Natural Language Processing}, pp. 1746–-1751.

\bibitem{KoKr15}
Kohler, M., and Krzy\.zak, A. (2017).
Nonparametric regression based on hierarchical interaction models.
\textit{IEEE Transaction on Information Theory} \textbf{63}: 1620-1630.




\bibitem{KoKr19}
  Kohler, M., and Krzy\.zak, A. (2021).
  Over-parametrized deep neural networks
  minimizing the empirical risk
  do not generalize well. 
  {\it Bernoulli} \textbf{27}(4): 2564-2597.



\bibitem{KoKrLa19}
  Kohler, M., Krzy\.zak, A., and Langer, S. (2022).
  Estimation of a function of low local dimensionality by deep neural
 networks.
\textit{IEEE Transactions on Information Theory} \textbf{68}(6): 4032 -- 4042.


  
  \bibitem{KoLa19}
  Kohler, M., and Langer, S. (2022).
  Discussion of ``Nonparametric regression using deep neural networks with ReLU activation function''.
  {\it Annals of Statistics} {\bf 48}(4): 1906-1910.
  

\bibitem{KoLa19}
  Kohler, M., and Langer, S. (2022).
  On the rate of convergence of fully connected very deep
  neural network regression estimates using ReLU activation functions.
 {\it Annals of Statistics} {\bf 49}(4): 2231--2249.



  
\bibitem{Kol04}
  Koltchinskii, V. (2006).
  Local Rademacher complexities and oracle inequalities
  in risk minimization. {\it Annals of Statistics}
  {\bf 34}(6):2593 -- 2656.




\bibitem{KSH12}
Krizhevsky, A., Sutskever, I., and Hinton, G. E. (2012).
ImageNet classification with deep convolutional neural networks.
In {\it Advances In Neural Information
Processing Systems} {\bf 25}, pp. 1097--1105. 

\bibitem{KS08}
K\r{u}rková,  V., and Sanguinetti, M. (2008). Geometric upper bounds on rates of variable-basis approximation. {\it IEEE Transactions on Information Theory} {\bf 54}: 5681 -- 5688.

\bibitem{Lee96}
  Lee, W.S. (1996). Agnostic Learning and
Single Hidden Layer Neural Networks. PhD Thesis, Australian National
University.


\bibitem{LeLiSrSu18}
Liang, S., Sun, R., Lee, J., and Srikant, R. (2018).
Adding one neuron can eliminate all bad local minima.
\textit{Advances In Neural Information Processing Systems},
pp. 4355 -- 4365.

 
\bibitem{LiLi18}
Li, Y., and Liang, Y. (2018).
Learning overparameterized neural networks via stochastic gradient descent on structured data. In
\textit{Advances in Neural Information Processing Systems},
pp. 8168- 8177.



\bibitem{McCaGa94}
McCaffrey, D. F., and Gallant, A. R. (1994).
Convergence rates for single hidden layer feedforward
networks.
{\it Neural Networks} {\bf 7}: 147-158.

\bibitem{Pisier}
Pisier, G.
"Remarques sur un resultat non publie de B. Maurey," presented at the Seminaire d'analyse fonctionelle 1980-1981, Ecole Polytechnique, Centre de Mathematiques, Palaiseau.

\bibitem{P}
Poggio, T., Banburski, A. , and Liao,Q.
Theoretical issues in deep networks.
{\it Proceedings of the National Academy of Sciences} {\bf 117}: 30039--30045.

\bibitem{Rahimi}
Rahimi, A., and Recht, B. (2009). 
Weighted sums of random kitchen sinks: Replacing minimization with randomization in learning.  In {\it Advances in Neural Information Process Systems} {\bf 21}, pp. 1313-1320. 

\bibitem{Sch17}
Schmidt-Hieber, J. (2020).
Nonparametric regression using deep neural networks with ReLU activation function (with discussion). {\it Annals of Statistics} {\bf 48}(4): 1875--1897. 


  
\bibitem{Setal17}
  Silver, D., Schrittwieser, J., Simonyan, K., Antonoglou, I., Huang, A.,
  Guez, A., Huber, T., et al. (2017).
  Mastering the game of go without human knowledge.
  {\it Nature} {\bf 550}: 354-359.

\bibitem{Sit20}
Sitzmann,  V., Martel, J., Bergman,A., Lindell, D., and Wetzstein, G. (2020).
Implicit neural representations with periodic activation functions. In
{\it Advances in Neural Information Processing Systems} {\bf 33}, pp. 7462--7473.

\bibitem{So21}
Sonoda, S., Ishikawa, I., and Ikeda, M. (2021). Ridge regression with overparametrized two-layer networks convergence to ridgelet spectrum. 
{\it International Conference on Artificial Intelligence and Statistics} {\bf 130}, pp. 2674 – 2682.

  
\bibitem{Sto82}
Stone, C. J. (1982).
Optimal global rates of convergence for nonparametric regression.
\textit{Annals of Statistics} \textbf{10}(4): 1040-1053.


\bibitem{SN19}
Suzuki, T., and Nitanda, A. (2021).
Deep learning is adaptive to intrinsic diemsnionality of
model smoothness in anisotropic Besov space.
{\it Advances in Neural Information Processing Systems} {\bf 34}, pp.  3609-3621.

\bibitem{Widrow}
Widrow, B., Greenblatt, A., Kim, Y., and Park, D. (2013). The No-Prop algorithm: A new learning algorithm for multilayer neural networks. {\it Neural Networks} {\bf 37}: 182-188.

\bibitem{Wetal20}
  Woodworth, B., Gunasekar, S., Lee, J., Moroshko, E., Savarese, P., Golan, I.,
Soudry, D.,  und  Srebro, N. (2020).   Kernel  and  rich  regimes  in  overparametrized  models. In
{\it Conference on Learning Theory} {\bf 125}, pp. 3635--3673.

  

\bibitem{Wetal16}
  Wu, Y., Schuster, M., Chen, Z., Le, Q., Norouzi, M.,
  Macherey, W., Krikum, M., et al. (2016).
  Google's neural machine translation system: Bridging the gap between
  human and machine translation. \textit{arXiv: 1609.08144}.
  

\bibitem{Yos68}
  Yosida, K. (1968). {\it Functional Analysis}. 2nd edition. Springer.
Berlin.
  

\bibitem{Zouetal18}
  Zou, D., Cao, Y., Zhou, D., und Gu, Q. (2020).
  Gradient descent optimizes over-parameterized
  deep ReLU networks.
 \textit{Machine Learning} {\bf 109}: 467 -- 492.

\end{thebibliography}
\end{document}